# MODULES OVER IWASAWA ALGEBRAS

J. Coates, P. Schneider, R. Sujatha

**Introduction**

Let $p$ be a prime number, and $G$ a compact $p$-adic Lie group. We recall that the Iwasawa algebra (or completed group ring of $G$) is defined by

$$\Lambda(G) = \varprojlim_{H} \mathbb{Z}_p[G/H],$$

where $H$ runs over the set of all open normal subgroups of $G$. Interesting examples of finitely generated modules over $\Lambda(G)$, in which $G$ is the image of Galois in the automorphism group of a $p$-adic Galois representation, abound in arithmetic geometry (see §7 for some examples coming from elliptic curves without complex multiplication). The study of such $\Lambda(G)$-modules arising from arithmetic geometry can be thought of as a natural generalization of Iwasawa theory. One of the cornerstones of classical Iwasawa theory is the fact that, when $G = \mathbb{Z}_p^d$ for some integer $d \geq 1$, a good structure theory for finitely generated $\Lambda(G)$-modules is known, up to pseudo-isomorphism. When $d = 1$, this was proven by Iwasawa [Iw] by an ad hoc method, but almost immediately Serre [Se0], [Se1] pointed out that it was a special case of one of the main results of commutative algebra and could be proved for all $d \geq 1$ (see [B-CA, Chap. VII, §4.4 Thms. 4 and 5]).

The aim of the present paper is to extend as much as possible of this commutative structure theory to the non-commutative case. The first important step in this direction was taken by Venjakob [Ven0], [Ven1], who, using ideas of Björk [Bjo], discovered an appropriate non-commutative generalization of the notion of a pseudo-null module over $\Lambda(G)$. Inspired by Venjakob's work, we propose in §1 a definition of pseudo-null modules over an arbitrary ring $A$; in particular, a finitely generated torsion $\Lambda(G)$-module $M$ is pseudo-null if $\mathrm{Ext}^1_{\Lambda(G)}(M, \Lambda(G)) = 0$ whenever $G$ is pro-$p$ and has no element of order $p$. Venjakob [Ven0], [Ven1] then went on to establish the first fragment of the structure theory by treating the case of finitely generated $\Lambda(G)$-modules which are annihilated by some power of $p$ (see also [Ho1], [Ho2]). In §3 of this paper, we exploit a well known filtration on $\Lambda(G)$ (see §6) to show that ideas about filtered rings, which have their origin in the algebraic theory of microlocalization and seem to go back to [Ka], do enable one to rather simply and elegantly extend Bourbaki's proof of the structure theory in the commutative case to finitely generated torsion modules over $\Lambda(G)$, for a wide class of $p$-adic Lie groups $G$. We must confess that we were greatly surprised at first to find that we could do this, as we had initially expected that the structure theory for non-abelian groups $G$ would be fundamentally different from the case of $G = \mathbb{Z}_p^d$.

We also then realized that this structure theory could be derived from earlier work of Chamarie [Ch0], [Ch1], [Ch2] on modules over maximal orders. We explain in detail the connection with Chamarie's work in §2 of this paper, but we do not repeat the proof of his version of the structure theory, which he formulates only somewhat less precisely than us in the quotient of the category of all $\Lambda(G)$-modules by the subcategory of pseudo-null submodules.

We assume that our compact $p$-adic Lie group $G$ is $p$-valued in the sense of Lazard [Laz] (we recall Lazard's definition in §6). We are very grateful to B. Totaro for pointing out to us that our arguments are valid in this generality, and that the class of $p$-valued compact $p$-adic Lie groups enjoys many nice properties (e.g. that every closed subgroup of such a group is again $p$-valued). Every compact $p$-adic Lie group contains an open subgroup which is $p$-valued. Important examples of $p$-valued $G$ are given by the subgroup of $GL_n(\mathbb{Z}_p)$ consisting of all matrices which are congruent to 1 modulo $p$ or modulo 4, according as $p$ is odd or even, and by any pro-$p$ closed subgroup of $GL_n(\mathbb{Z}_p)$ provided $p > n+1$. The precise statement of our Structure Theorem is as follows:

**Structure Theorem:** *Let $G$ be a $p$-valued compact $p$-adic Lie group, and let $M$ be a finitely generated torsion $\Lambda(G)$-module. Let $M_\mathrm{o}$ be the maximal pseudo-null submodule of $M$. Then there exist non-zero left ideals $L_1, \cdots, L_m$, and a $\Lambda(G)$-injection*

$$\phi \;:\; \bigoplus_{i=1}^m \Lambda(G)/L_i \;\longrightarrow\; M/M_\mathrm{o}$$

*with* $\mathrm{Coker}(\phi)$ *pseudo-null.*

We strongly believe there is genuine interest in retaining a rather full discussion of these two somewhat different approaches to the Structure Theorem in the present paper. We exploit Chamarie's methods in §4 to prove some form of uniqueness for the ideals $L_1, ..., L_m$ appearing in the Structure Theorem, and to define the notion of the characteristic ideal of $M$. We also study arbitrary finitely generated $\Lambda(G)$-modules in §5 by Chamarie's methods. We have also thought it instructive to present our proofs throughout in the axiomatic fashion, leaving the verification that $\Lambda(G)$ satisfies these axioms until §6. Finally, we point out that perhaps the major question left open by our work is whether the left ideals $L_1, \cdots, L_m$ which occur in the above structure theorem can be chosen to be principal (when $G = \mathbb{Z}_p^d$ this is true because $\Lambda(G)$ is a unique factorization domain).

In conclusion, we would like to express our warmest thanks both to O. Venjakob and S. Howson for keeping us fully informed on their work. The first and the third author are also deeply grateful to H. Li and F. van Oystaeyen for their infinite patience in answering their inexpert questions on the algebraic theory of microlocalization. The first author wishes to thank both Columbia University and the long suffering audience of his Eilenberg lectures there in the spring of 2001 for supporting this research so generously. The second author gratefully acknowledges the support of the Tata Institute.

**Notation:** We summarize here the most important notation used in the paper. All of our rings will be assumed associative and with unit element, but will not, in general be commutative. If $A$ is such a ring, we will always consider left $A$-modules unless the contrary is stated. We write $\text{Mod}(A)$ for the category of left $A$-modules, and $\mathfrak{M}(A)$ for the subcategory of finitely generated left $A$-modules. A module $M$ in $\text{Mod}(A)$ will be said to be *A-torsion* if every element of $M$ has a non-zero annihilator in $A$. We shall often assume that we are given a filtration on $A$, written $F.A = \{F_n A : n \in \mathbb{Z}\}$, which we shall always assume is indexed by $\mathbb{Z}$, increasing, and exhaustive. We write $\text{gr}.A = \oplus_{n \in \mathbb{Z}} F_n A / F_{n-1} A$ for the associated graded ring. We use a similar notation and convention for filtrations on left $A$-modules. We say that $A$ is noetherian if every left ideal and right ideal in $A$ is finitely generated. If $K$ is any field, we write $\bar{K}$ for a separable closure of $K$. If $p$ is a prime number, $\mathbb{Z}_p$ will denote the ring of $p$-adic integers, and $\mathbb{Q}_p$ the quotient field of $\mathbb{Z}_p$.

## 1. Pseudo-null modules

As above, let $A$ be an associative ring with unit and let $\text{Mod}(A)$ be the category of left $A$-modules. For the moment, we impose no further conditions on $A$. One way to introduce in this great generality a dimension filtration on the category $\text{Mod}(A)$ is through the following technique. As usual $E(M)$ denotes the injective hull (compare [Ste, Chap. V, §2]) of the $A$-module $M$. Consider the "minimal" injective resolution

$$0 \longrightarrow L \xrightarrow{\mu_0} E_0 \xrightarrow{\mu_1} E_1 \xrightarrow{\mu_2} \ldots$$

of the $A$-module $L$, i.e., such that $E_0 = E(L)$ and $E_{i+1} = E(\text{Coker}(\mu_i))$ for any $i \geq 0$. Let

$$\mathcal{C}_L^n := \text{ full subcategory of all } M \text{ in } \text{Mod}(A)$$
$$\text{such that } \text{Hom}_A(M, E_0 \oplus \ldots \oplus E_n) = 0 \ .$$

This subcategory $\mathcal{C}_L^n$ is "localizing" in the sense that it satisfies the following conditions:

– In any short exact sequence $0 \to M' \to M \to M'' \to 0$ of $A$-modules $M$ lies in $\mathcal{C}_L^n$ if and only $M'$ and $M''$ lie in $\mathcal{C}_L^n$;

– any $A$-module has a unique largest submodule contained in $\mathcal{C}_L^n$.

It is called the *hereditary torsion theory* cogenerated by the injective module $E_0 \oplus \ldots \oplus E_n$ (see [Ste, Chap. VI]).

**Lemma 1.1.** *An $A$-module $M$ lies in $\mathcal{C}_L^n$ if and only if $\operatorname{Ext}_A^i(M', L) = 0$ for any $i \leq n$ and any submodule $M' \subseteq M$.*

*Proof.* See [Ste, Chap. VI, Prop. 6.9]. □

Suppose that $A$ is noetherian, and has no divisors of zero. Then the classical (left and right) ring of quotients $D$ of $A$ exists and is a skew-field [Ste, Chap. II, Prop. 1.7]. Since $D$ is injective as a left and right $A$-module [Ste, Chap. II, Prop. 3.8] we have $E(A) = D$. Further, as

$$\operatorname{Hom}_A(M, D) = \operatorname{Hom}_D(D \otimes_A M, D) = 0 \quad \text{if and only if} \quad D \otimes_A M = 0$$

it follows that

$$\begin{aligned} \mathcal{C}_A^0 \;\; &= \;\; \text{full subcategory of all torsion modules } M \\ & \quad \text{(i.e., such that } D \otimes_A M = 0\text{)} \; . \end{aligned}$$

**Lemma 1.2.** *Suppose that $A$ is noetherian without zero divisors and that $A \neq D$; we then have*

$$\begin{aligned} \mathcal{C}_A^1 \;\; &= \;\; \text{full subcategory of all } M \text{ such that} \\ & \quad \operatorname{Hom}_A(M', D/A) = 0 \text{ for any submodule } M' \subseteq M \; . \end{aligned}$$

*Proof.* Since $D/A$ is an essential submodule of $E_1 = E(D/A)$ the asserted condition on $M$ is equivalent to $\operatorname{Hom}_A(M, E_1) = 0$. It therefore remains to show that this condition implies that $\operatorname{Hom}_A(M, D) = 0$. It certainly implies that

$$\operatorname{Hom}_A(M, A) = \operatorname{Hom}_A(M, D) \; .$$

Consider now any $f \in \operatorname{Hom}_A(M, A)$. For any non-zero $b \in A$ we have the $A$-linear map $f_b : M \longrightarrow D$ defined by $f_b(x) := f(x)b^{-1}$. This shows that

$$f(M) \subseteq \bigcap_{0 \neq b \in A} Ab \; .$$

Suppose that the right hand side contains an $a \neq 0$. Choose a non-zero and non-invertible $c \in A$ ($A \neq D$!) and set $b := ca \neq 0$. There must exist a $d \in A$ such that $a = db = dca$ and hence $1 = dc$ which is a contradiction. □

**Example 1.3.** (cf. [Ste, Chap. VII, Prop. 6.10]) Suppose that $A$ is an integrally closed noetherian commutative integral domain. Then an $A$-module $M$ lies in $\mathcal{C}_A^1$ if and only if $M_\mathfrak{p} = 0$ for any prime ideal $\mathfrak{p} \subseteq A$ of height $\leq 1$. Hence a finitely generated $A$-module lies in $\mathcal{C}_A^1$ if and only if it is pseudo-zero in the sense of [B-CA, Chap. VII, §4.4].

Another very interesting class of rings (which includes Iwasawa algebras, see [Ven0], [Ven1]) consists of the Auslander regular rings, and we now recall their definition. The *grade* of a left or right $A$-module $M$ is defined to be the smallest non-negative integer $j(M) = j_A(M)$ such that $\operatorname{Ext}_A^{j(M)}(M, A) \neq 0$ (we let $j(\{0\}) = \infty$). For any finitely generated $M \neq 0$ the grade $j(M)$ is bounded above by the projective dimension of $M$. We say that $M$ satisfies the *Auslander condition* if, for each $k \geq 0$ and any submodule $N$ of $\operatorname{Ext}_A^k(M, A)$, we have $j_A(N) \geq k$. A noetherian ring $A$ is defined to be *Auslander regular* if every finitely generated left or right $A$-module has finite projective dimension and satisfies the Auslander condition (this is slightly more general than that given in [LVO, Chap. III, Def. 2.1.7] in that we do not necessarily assume that $A$ has finite global dimension). We mention that a commutative noetherian ring is Auslander regular if and only if it is regular (see [AB, Cor. 4.6, Prop. 4.21] or [BH, Cor. 3.5.11]). Returning to a general Auslander regular ring $A$ we note that any finitely generated $A$-module $M$ lies in the subcategory $\mathcal{C}_A^{j(M)-1}$. We say that a module $M$ is *pure* if $\operatorname{Ext}_A^i(\operatorname{Ext}_A^i(M, A), A) = 0$ for any $i \neq j(M)$. Suppose that the $A$-module $M$ is finitely generated and let $d$ denote its projective dimension. Then $M$ (see [Bjo] or [LVO, Chap. III]) carries a natural filtration, called the *dimension filtration*, by submodules

$$M = \Delta^0(M) \supseteq \Delta^1(M) \supseteq \ldots \supseteq \Delta^{d+1}(M) = 0$$

(note that we have changed the numbering of this filtration compared to [Bjo] and [LVO], because we believe that the above, in which the numbering corresponds to codimension, is more natural). This filtration is characterized by the property that a submodule $L \subseteq M$ has grade $j(L) \geq p$ if and only if $L \subseteq \Delta^p(M)$. In addition one has:

– $j(M) = \max\{p \geq 0 : \Delta^p(M) = M\}$ ;

– if $M$ is pure then $M = \Delta^{j(M)}(M) \supset \Delta^{j(M)+1}(M) = 0$ ;

– $\Delta^p(M)/\Delta^{p+1}(M)$ is zero or pure of grade $p$ .

**Lemma 1.4.** *If $A$ is Auslander regular then a finitely generated $A$-module $M$ lies in $\mathcal{C}_A^n$ if and only if $j_A(M) > n$.*

*Proof.* Suppose that $j(M) > n$. We know that $j(M') \geq j(M) > n$ for any submodule $M' \subseteq M$ which by Lemma 1.1 implies that $M$ lies in $\mathcal{C}_A^n$. The converse is immediate from Lemma 1.1. □

This result shows in particular that over an Auslander regular ring $A$ a finitely generated module $M$ lies in $\mathcal{C}_A^1$ if and only it is pseudo-null in the sense of [Ven0, 1.5.1]. In the light of these two examples we suggest the following definition.

**Definition 1.5.** A module $M$ over an arbitrary ring $A$ is called *pseudo-null* if it lies in $\mathcal{C}_A^1$.

## 2. Review of Chamarie's work

We assume throughout this section that $A$ is noetherian and has no zero divisors, and we write $D$ for the skew-field of fractions of $A$. If $A$ is commutative and integrally closed, then one of the basic results of commutative algebra is the structure theory of finitely generated torsion $A$-modules up to pseudo-isomorphism [B-CA, Chap. VII, §4.4 Thm. 5]. In this section, we discuss one generalization of this theory to non-commutative rings due to Chamarie in [Ch2], and in the next section, we shall present in detail a second generalization based on ideas from the algebraic theory of microlocalization.

We recall that $A$ is called a *maximal order* if, for any intermediate ring $B$ with $A \subseteq B \subseteq D$ such that there exist elements $u, v \in D^\times$ with $uBv \subseteq A$, we always have $A = B$.

**Examples 2.1.** (cf. [MCR, Chap. 5, §1 and §3])

∗ If $A$ is commutative then it is a maximal order if and only if it is integrally closed.

∗ Any Weyl algebra over a field $k$ and any universal enveloping algebra of a finite dimensional Lie algebra over $k$ is a maximal order.

∗ Let $R$ be a commutative noetherian integral domain with quotient field $k$; any $R$-order in a central simple algebra over $k$ which is maximal with respect to inclusion is a maximal order (possibly with zero divisors).

For later use we quote the following criterion.

**Lemma 2.2.** *The ring $A$ is a maximal order if and only if for any non-zero two sided ideal $I \subseteq A$ and any $u \in D \setminus A$ we have $uI \nsubseteq I$ and $Iu \nsubseteq I$.*

*Proof.* See [MCR, Prop. 5.1.4]. □

We write $\mathrm{Mod}(A)/\mathcal{C}_A^1$ for the quotient category of $\mathrm{Mod}(A)$ by the subcategory $\mathcal{C}_A^1$ of pseudo-null $A$-modules, and put

$$Q : \mathrm{Mod}(A) \longrightarrow \mathrm{Mod}(A)/\mathcal{C}_A^1$$

for the quotient functor. We call an $A$-module $M$ *pseudo-cyclic* if it contains a cyclic submodule $N$ such that the quotient $M/N$ lies in $\mathcal{C}_A^1$. It turns out to be rather involved to define the notion of a cyclic object in the quotient category $\mathrm{Mod}(A)/\mathcal{C}_A^1$. We define an object $Q(M)$ in $\mathrm{Mod}(A)/\mathcal{C}_A^1$ to be *cyclic* if every $A$-module $N$ such that $Q(N) \cong Q(M)$ is pseudo-cyclic. According to [Rob, Lemma 1.6], the object $Q(M)$ is cyclic if and only if any submodule $N \subseteq M$ such that $M/N$ lies in $\mathcal{C}_A^1$ contains a cyclic submodule $L$ such that $M/L$ lies in $\mathcal{C}_A^1$. In general it need not even be true that $Q(A)$ is cyclic, but, for the rings which interest us, we shall show later (Lemma 2.7) that the notion of $Q(M)$ being cyclic for an $A$-torsion module $M$ coincides with what we intuitively hope to be true. In the quotient category $\mathrm{Mod}(A)/\mathcal{C}_A^1$ Chamarie establishes the following important structure theorem when $A$ is a maximal order.

**Proposition 2.3.** *If $A$ is a maximal order then any object of finite length in $\mathrm{Mod}(A)/\mathcal{C}_A^1$ is a direct sum of cyclic objects.*

*Proof.* See [Ch2, Thm. 4.2.7]. □

We next explain how one can describe the subobjects of $Q(A)$ in terms of a certain class of left ideals of $A$. For an arbitrary hereditary torsion theory $\mathcal{C}$ in $\mathrm{Mod}(A)$, one has the notion of a $\mathcal{C}$-closed left ideal in $A$ [Ste, Chap. IX, §1]. We do not recall the definition here, since we can work with the following characterization (cf. [Ste, Chap. IX, Prop. 4.2]), because the torsion theory $\mathcal{C}_A^1$ has the property that the ring $A$ itself is $\mathcal{C}_A^1$-closed [Ste, Chap. IX, Prop. 2.15]. A left ideal $L \subseteq A$ is $\mathcal{C}_A^1$-*closed* if and only if the quotient module $A/L$ is $\mathcal{C}_A^1$-torsion free, i.e., does not contain any non-zero submodule lying in $\mathcal{C}_A^1$. Granted this characterization, we see immediately that, for the torsion theory $\mathcal{C}_A^1$, the definitions of a $\mathcal{C}_A^1$-closed left ideal given in [Ste] and in [Ch0-2] coincide.

If $A$ is Auslander regular, the properties of the dimension filtration show that a finitely generated $A$-module is $\mathcal{C}_A^1$-torsion free if and only if $\Delta^2(M) = 0$. In particular, a left ideal $L \subseteq A$ is $\mathcal{C}_A^1$-closed if and only if $\Delta^2(A/L) = 0$. Since for non-zero $L$ the quotient $A/L$ is $A$-torsion it follows that a proper left ideal $0 \subset L \subset A$ is $\mathcal{C}_A^1$-closed if and only if the quotient module $A/L$ is pure of grade 1.

The crucial observation in describing the subobjects of $Q(A)$ is the following [Ste, Chap. IX, Cor. 4.4], in which we have once more dropped the assumption that $A$ is Auslander regular. The quotient functor $Q$ induces a bijection between the set of $\mathcal{C}_A^1$-closed left ideals in $A$ and the set of subobjects of $Q(A)$. We denote by $K_1\dim(A)$ the Krull dimension of the object $Q(A)$ in the quotient category $\mathrm{Mod}(A)/\mathcal{C}_A^1$. Equivalently, $K_1\dim(A)$ is the deviation [MCR, Chap. 6, §1] of the partially ordered (w.r.t. inclusion) set of $\mathcal{C}_A^1$-closed left ideals in $A$.

We shall be mainly interested in rings $A$ with the property that $K_1\dim(A) \leq 1$. By definition this means that for any descending chain $L_1 \supseteq L_2 \supseteq \ldots$ of $\mathcal{C}_A^1$-closed left ideals $L_i \subseteq A$ all but finitely many of the sets

$$\{L' \subseteq A : L' \text{ a } \mathcal{C}_A^1\text{-closed left ideal with } L_{i+1} \subseteq L' \subseteq L_i\}$$

satisfy the descending chain condition. But actually a more precise statement holds true. Let $L \subseteq A$ be a non-zero $\mathcal{C}_A^1$-closed left ideal. Suppose that

$$A \supset L_1 \supset L_2 \supset \ldots \supset L$$

is a strictly descending infinite chain of $\mathcal{C}_A^1$-closed left ideals. We fix a non-zero element $b \in L$ and consider the descending chain

$$Q(A) \supseteq Q(Ab) \supseteq Q(Ab^2) \supseteq \ldots$$

of subobjects of $Q(A)$ (recall that the quotient functor $Q$ is exact). For each $j \geq 0$ we have the intermediate chain

$$Q(Ab^j) \supset Q(L_1 b^j) \supset Q(L_2 b^j) \supset \ldots \supseteq Q(Lb^j) \supseteq Q(Ab^{j+1}) .$$

Since

$$Q(L_i b^j)/Q(L_{i+1} b^j) = Q(L_i b^j / L_{i+1} b^j) \cong Q(L_i/L_{i+1}) \neq 0$$

all inclusions in these latter chains labeled "⊃" are strict. This shows that $K_1\dim(A) > 1$. The condition $K_1\dim(A) \leq 1$ therefore ensures that the object $Q(A/L)$, for any non-zero $\mathcal{C}_A^1$-closed left ideal $L \subseteq A$, has finite length. By a simple induction on the number of generators we derive the following fact.

**Remark 2.4.** *Suppose that $K_1\dim(A) \leq 1$; for any finitely generated torsion $A$-module $M$ the object $Q(M)$ has finite length in $\mathcal{C}_A^0/\mathcal{C}_A^1$.*

**Corollary 2.5.** *Suppose that $A$ is a maximal order satisfying $K_1\dim(A) \leq 1$ and let $M$ be any finitely generated torsion $A$-module. We then have:*
*(i) $Q(M)$ is a direct sum of cyclic objects.*
*(ii) There are non-zero $\mathcal{C}_A^1$-closed left ideals $L_1, \ldots, L_m \subseteq A$ together with an injective homomorphism of $A$-modules with pseudo-null cokernel*

$$\bigoplus_{i=1}^{m} A/L_i \longrightarrow M/M_\mathrm{o}$$

*where $M_\mathrm{o}$ denotes the largest pseudo-null submodule of $M$.*

*Proof.* *(i):* This is immediate from Prop. 2.3 and the above Remark 2.4.

*(ii):* By *(i)* we have an isomorphism

$$\bigoplus_{i=1}^{m} Q(N_i) \xrightarrow{\cong} Q(M)$$

in the quotient category $\mathrm{Mod}(A)/\mathcal{C}_A^1$ where the $Q(N_i)$ are cyclic objects. Consider the restriction of this morphism to a fixed summand on the left hand side

$$Q(N_i) \longrightarrow Q(M) \ .$$

By the construction of the quotient category it is represented by an actual homomorphism of $A$-modules

$$N_i' \longrightarrow M/M_i'$$

where $N_i' \subseteq N_i$ and $M_i' \subseteq M$ are submodules such that $N_i/N_i'$ and $M_i'$ are pseudo-null. In particular, we have $M_i' \subseteq M_\mathrm{o}$. Since $N_i$ is pseudo-cyclic we may assume that $N_i'$ is a cyclic $A$-module, i.e., is of the form $N_i' = A/L_i'$ for some non-zero left ideal $L_i' \subseteq A$. Let $L_i \subseteq A$ be the left ideal such that $L_i/L_i'$ is the largest pseudo-null submodule of $A/L_i'$. Then $L_i$ is $\mathcal{C}_A^1$-closed and the above map induces a homomorphism

$$A/L_i \longrightarrow M/M_\mathrm{o} \ .$$

The corresponding sum

$$\bigoplus_{i=1}^{m} A/L_i \longrightarrow M/M_\mathrm{o}$$

represents in the quotient category the original isomorphism and therefore has pseudo-null kernel and cokernel. Since the left hand side, by construction, has no non-zero pseudo-null submodule this last map actually is injective. □

According to [Ch2, Cor. 4.1.2] the condition $\mathrm{K}_1\mathrm{dim}(A) \leq 1$ is satisfied if $A$ is bounded (i.e., every non-zero left or right ideal contains a non-zero two sided ideal) and hence in particular if $A$ is commutative.

To provide further examples of maximal orders including the concrete examples of Iwasawa algebras, and which satisfy $\mathrm{K}_1\mathrm{dim}(A) \leq 1$, we now consider a class of filtered rings $A$. Thus we now assume that $A$ is endowed with an exhaustive and complete filtration $F.A$ (we observe the conventions on filtrations as explained in the Notations) such that the associated graded ring $\mathrm{gr}.A$ is noetherian without zero divisors. The reader should compare this assumption with the axioms (C1)–(C3) imposed on $A$ in §3, noting in particular that we not impose here that gr.$A$ is commutative. As is explained in §3 (see also §3 for an explanation of a good filtration on an $A$-module), the above assumption guarantees that $A$ is a Zariski ring in the sense of [LVO]. In particular, $A$ has the following properties:

– $A$ is noetherian without zero divisors [LVO, Chap. II, Prop. 1.2.3].

– Any left ideal $L \subseteq A$ is closed in the filtration topology [LVO, Chap. I, Cor. 5.5].

– The filtration induced by $F.A$ on any subquotient module of $A$ is good [LVO, Chap. II Prop. 1.2.3].

– If gr.$A$ is Auslander regular, then so is the ring $A$ itself [LVO, Chap. III, Thm. 2.2.5].

– If $M$ is any finitely generated $A$-module endowed with a good filtration, we then have [LVO, Chap. III, Thm. 2.5.2]

$$j_A(M) = j_{\text{gr}.A}(\text{gr}.M) \ .$$

**Lemma 2.6.** *Assume that $A$ is endowed with an exhaustive and complete filtration such that* gr.$A$ *is noetherian without zero divisors. Then:*
*(i) If* gr.$A$ *is a maximal order, so is $A$.*
*(ii) If* gr.$A$ *is Auslander regular, then* $\text{K}_1\text{dim}(A) \leq \text{K}_1\text{dim}(\text{gr}.A)$.

*Proof.* (i): This is [Ch1, Prop. 2.3.1]. For the convenience of the reader we recall the argument. We will use the criterion of Lemma 2.2. Let $I \subseteq A$ be a non-zero two sided ideal and let $u \in D$ be such that $uI \subseteq I$. If $u = b^{-1}a$ with $a, b \in A$ then the latter is equivalent to $aI \subseteq bI$. We have to show that $u \in A$, i.e., that $a \in bA$. That the ideal $bA$ is closed in the filtration topology means that

$$bA = \bigcap_{n \in \mathbb{Z}} bA + F_n A \ .$$

Hence it suffices to show that $a \in bA + F_n A$ for all sufficiently small $n \in \mathbb{Z}$. The filtration being exhaustive we find an $m \in \mathbb{Z}$ such that $a \in F_m A \subseteq bA + F_m A$. So it furthermore suffices to prove that $a \in bA + F_{n-1} A$ provided $a \in bA + F_n A$. Let us assume that $a = by + z$ with $y \in A$ and $z \in F_n A$. Obviously we only need to treat the case where $z \neq 0$. Suppose that $z \in F_m A \setminus F_{m-1} A$ and $b \in F_p A \setminus F_{p-1} A$ and put $\overline{z} := z + F_{m-1} A \in \text{gr}.A$ and $\overline{b} := b + F_{p-1} A \in \text{gr}.A$. We have

$$zI \subseteq aI + byI \subseteq bI \ .$$

Forming gr.$I$ with respect to the induced filtration gives a non-zero two sided ideal in gr.$A$. We claim that

$$\overline{z} \cdot \text{gr.}I \subseteq \overline{b} \cdot \text{gr.}I \ .$$

Consider any $0 \neq \overline{c} = c + F_{r-1}A \in \text{gr}_r I$. Then $zc = bd$ for some $0 \neq d \in I$. Suppose that $d \in F_s I \setminus F_{s-1} I$ and put $\overline{d} := d + F_{s-1}I \in \text{gr}_s I$. Since gr.$A$ has no zero divisors we must have $m + r = p + s$ and

$$\overline{z} \cdot \overline{c} = zc + F_{m+r-1}A = bd + F_{p+s-1}A = \overline{b} \cdot \overline{d} \ .$$

This proves our claim. Since gr.$A$ is a maximal order it follows that

$$\overline{z} \in \overline{b} \cdot \text{gr.}A \ .$$

Hence there is a $\overline{x} = x + F_{s-r-1}A \in \text{gr}_{s-r-1}A$ such that $\overline{z} = \overline{b} \cdot \overline{x} = bx + F_{m-1}A$. We obtain $z - bx \in F_{m-1}A \subseteq F_{n-1}A$ and consequently

$$a = by + z = b(y + x) + z - bx \in bA + F_{n-1}A \ .$$

The argument for $Iu \subseteq I$ implying that $u \in A$ is analogous.

(ii): Let $0 \subset L \subset A$ be a proper $\mathcal{C}_A^1$-closed left ideal. Equipping $L$ and $A/L$ with the induced filtrations we obtain the proper left ideal gr.$L \subset$ gr.$A$ such that $j(\text{gr.}A/\text{gr.}L) = j(\text{gr.}A/L) = j(A/L) = 1$. Let gr.$L \subset \widetilde{\text{gr}}L \subset$ gr.$A$ denote the proper left ideal such that $\widetilde{\text{gr}}L/\text{gr.}L = \Delta^2(\text{gr.}A/\text{gr.}L)$. Then gr.$A/\widetilde{\text{gr}}L$ is pure of grade 1 so that $\widetilde{\text{gr}}L$ is $\mathcal{C}_{\text{gr.}A}^1$-closed. Hence $L \longmapsto \widetilde{\text{gr}}L$ is an inclusion preserving map from the set of proper $\mathcal{C}_A^1$-closed left ideals in $A$ to the set of proper $\mathcal{C}_{\text{gr.}A}^1$-closed left ideals in gr.$A$. To establish our assertion it suffices to show that this map preserves strict inclusions. Consider therefore two proper $\mathcal{C}_A^1$-closed left ideals $L \subset L'$ strictly contained in each other. By the purity of $A/L$ we have $j(\text{gr.}L'/\text{gr.}L) = j(\text{gr.}L'/L) = j(L'/L) = 1$. Hence gr.$L' \not\subseteq \widetilde{\text{gr}}L$ and a fortiori $\widetilde{\text{gr}}L' \neq \widetilde{\text{gr}}L$. □

This lemma applies, for example, if gr.$A$ is a commutative noetherian regular integral domain. Then gr.$A$ is integrally closed and hence a maximal order and is Auslander regular [LVO, Chap. III, 2.4.3] with $\text{K}_1\dim(\text{gr.}A) \leq 1$. Hence $A$ is an Auslander regular maximal order with $\text{K}_1\dim(A) \leq 1$.

Finally we note the following simple criterion for cyclicity in the quotient category $\text{Mod}(A)/\mathcal{C}_A^1$.

**Lemma 2.7.** *Suppose that $A$ is a maximal order satisfying $\text{K}_1\dim(A) \leq 1$ and let $M$ be a torsion $A$-module. Then $Q(M)$ is cyclic in $\text{Mod}(A)/\mathcal{C}_A^1$ if and only if there is a non-zero left ideal $L \subseteq A$ such that $Q(M) \cong Q(A/L)$.*

*Proof.* If $Q(M)$ is cyclic then, by definition, we find a cyclic torsion $A$-module $N$ such that $Q(M) \cong Q(N)$. But any such $N$ is isomorphic to $A/L$ for some non-zero left ideal $L \subseteq A$. For the reverse implication we use that, according to [Ch2, Thm. 4.2.4], any quotient of $Q(A)$ which is of finite length is cyclic. Under our assumptions $Q(A/L)$ is such a quotient by Remark 2.4. □

## 3. The approach via filtered rings

The aim of this section is to show that the classical proof of the structure theory of finitely generated torsion modules over commutative noetherian integrally closed domains (see [B-CA]) has a natural extension to a wide class of non-commutative filtered rings $A$. Thus, throughout this section $A$ will denote a ring, which is endowed with a separated and exhaustive filtration $F.A$, and which satisfies the following axioms, where, as usual, gr.$A$ denotes the associated graded ring:
(C1) $A$ is complete with respect to the filtration $F.A$;
(C2) gr.$A$ is a commutative noetherian integral domain;
(C3) gr.$A$ is a factorial ring in the sense of [B-CA, Chap. VII, §3], (i.e. gr.$A$ is integrally closed and every divisorial ideal is principal).

Our proof will make extensive use of techniques which grew out of the algebraic theory of microlocalization, and are, for example, studied systematically in [LVO].

Let $R$ be any filtered ring, with increasing filtration $F.R$. We recall (see [LVO, Chap. II, Thm. 2.2]) that one of the equivalent definitions of $R$ being a *Zariski ring* is that gr.$R$ is left and right noetherian, and that the completion $\hat{R} = \varprojlim_n R/F_n R$ is a faithful and flat left and right $R$-module. Applying this definition to our ring $A$, the next lemma is plain from axioms (C1) and (C2).

**Lemma 3.1.** *$A$ is a Zariski ring.* □

Throughout this section we fix a torsion module $M$ in $\mathfrak{M}(A)$. We shall always use a *good* filtration on $M$, whose definition we now recall. Pick any finite set $w_1, \cdots, w_s$ of generators of $M$ as an $A$-module, and fixed integers $k_1, \cdots, k_s$. For each $n$ in $\mathbb{Z}$, we then define

$$F_n M = \sum_{i=1}^{s} F_{n-k_i} A \cdot w_i.$$

It is clear that the $F_n M$ are an increasing sequence of additive subgroups of $M$, whose union is $M$. Any filtration of $M$ obtained in this manner is defined to be a good filtration. Let $F.M$ be a good filtration on $M$. As the ring $A$ is Zariski, it follows from the basic properties of these rings [LVO, Chap. I, Cor. 5.5 and Chap. II, Thm. 2.2] that not only is $F.M$ separated, but more generally that it satisfies the closure property, i.e. any submodule $N$ is closed in the filtration topology. Since $M$ is $A$-torsion, it is plain that gr.$M$ is gr.$A$-torsion. Also, gr.$M$ is a finitely generated gr.$A$-submodule because the filtration $F.M$ is good [LVO, Chap. I, Lemma 5.4].

By axiom (C2), gr.$A$ is a commutative ring, and our starting point for the proof of the structure theorem is to apply the classical commutative theory to the finitely generated torsion gr.$A$-module gr.$M$. For this commutative theory, we follow, wherever possible, the notation and terminology of [B-CA]. Hence we write Ass(gr.$M$) for the set of prime ideals $\mathfrak{p}$ of gr.$A$ which are associated to gr.$M$; recall that this set consists of prime ideals in gr.$A$ which are of the form ann($x$) for a non-zero element $x \in$ gr.$M$. Note that the zero ideal is not in Ass(gr.$M$) as gr.$M$ is a finitely generated torsion module. We define

$$W(M) = \{\mathfrak{p}_1, \cdots, \mathfrak{p}_r\}$$

to be the set of prime ideals of height one in Ass(gr.$M$).

**Remark 3.2.** While Ass(gr.$M$) in general depends on the choice of the good filtration $F.M$, it is true, as our notation suggests, that $W(M)$ is independent of this choice. More generally, let $U(M)$ denote the set of minimal elements of Ass(gr.$M$), so that $W(M)$ is a subset of $U(M)$. Then we claim that $U(M)$ is independent of the choice of the good filtration $F.M$. Indeed, let Supp(gr.$M$) denote the support of gr.$M$, i.e. the prime ideals of gr.$A$ which contain the annihilator ideal $\mathrm{ann}_{\mathrm{gr}.A}(\mathrm{gr}.M)$. Then it is well-known [B-AC, Chap. IV, §1.4, Thm. 2] that $U(M)$ coincides with the set of minimal elements of Supp(gr.$M$). Let $J(M)$ denote the radical of $\mathrm{ann}_{\mathrm{gr}.A}(\mathrm{gr}.M)$. This ideal is well-defined (cf. [Bjo, Prop. 5.1]) and independent of the filtration. As Supp(gr.$M$) is precisely the set of all prime ideals containing $J(M)$, the set $U(M)$ is independent of the filtration by the remark above, and hence so is $W(M)$.

**Definition 3.3.** The multiplicatively closed set $S$ is defined to be the set of all homogeneous elements in gr.$A$ which do not belong to $\mathfrak{p}_1 \cup \cdots \cup \mathfrak{p}_r$.

If $W(M)$ is empty, then $S$ is defined to be the set of all non-zero homogeneous elements of gr.$A$. As $S$ is a homogeneous multiplicatively closed set, we can clearly localize graded modules with respect to $S$. We also have the graded domain [B-CA, Chap. II, §2.9]

$$(\mathrm{gr}.A)_S := S^{-1}\mathrm{gr}.A.$$

Since gr.$M$ is a graded gr.$A$-module, every ideal in Ass(gr.$M$) is graded [B-CA, Chap. IV, §3.1, Prop. 1]. In particular, the ideals in $W(M)$ are all graded.

**Proposition 3.4.** *(i) The ideals $S^{-1}\mathfrak{p}_i$, $1 \leq i \leq r$, are graded prime ideals of height 1 in $(\mathrm{gr.}A)_S$.*
*(ii) Every proper graded ideal of $(\mathrm{gr.}A)_S$ is contained in one of the $S^{-1}\mathfrak{p}_i$, $1 \leq i \leq r$.*
*(iii) Every non-zero graded prime ideal of $(\mathrm{gr.}A)_S$ is of the form $S^{-1}\mathfrak{p}_i$ for some $i$, $1 \leq i \leq r$.*

*Proof. (i):* To establish *(i)*, we use [B-CA, Chap. II, §2.5, Prop. 11]. Now, for each $i = 1, \cdots, r$, $\mathfrak{p}_i$ is a prime ideal of $\mathrm{gr.}A$, which does not intersect $S$. Hence $S^{-1}\mathfrak{p}_i$ is a prime ideal of $(\mathrm{gr.}A)_S$. Moreover, if $\mathfrak{q}$ were any non-zero prime ideal of $(\mathrm{gr.}A)_S$ with $\mathfrak{q} \subseteq S^{-1}\mathfrak{p}_i$, we would necessarily have $\mathfrak{q} = S^{-1}\mathfrak{p}$ with $\mathfrak{p} \subseteq \mathfrak{p}_i$, and thus $\mathfrak{p} = \mathfrak{p}_i$ because $\mathfrak{p}_i$ has height one.

*(ii):* Turning to the proof of *(ii)*, let $\mathfrak{b}$ be any non-zero proper graded ideal of $(\mathrm{gr.}A)_S$. The preimage $\mathfrak{a}$ of $\mathfrak{b}$ in $\mathrm{gr.}A$ is a graded ideal such that $\mathfrak{b} = S^{-1}\mathfrak{a}$. Let $x_1, \cdots, x_k$ be a set of homogeneous generators for $\mathfrak{a}$ as a $\mathrm{gr.}A$-module, so that $x_1, \cdots, x_k$ also generate $\mathfrak{b}$ as a $(\mathrm{gr.}A)_S$-module. Since $\mathfrak{b}$ is a proper ideal of $(\mathrm{gr.}A)_S$, we must have $\mathfrak{a} \cap S = \emptyset$. Hence, as each $x_j$ is homogeneous, it must belong to at least one of the prime ideals $\mathfrak{p}_i$, and so we see that

$$\mathfrak{a} \subseteq \mathfrak{p}_1 \cup \cdots \cup \mathfrak{p}_r .$$

But, as each of the $\mathfrak{p}_i$ for $1 \leq i \leq r$ is a prime ideal, it follows from [B-CA, Chap. II, §1.1, Prop. 2] that $\mathfrak{a} \subseteq \mathfrak{p}_i$ for some $i$. Hence $\mathfrak{b} \subseteq S^{-1}\mathfrak{p}_i$, establishing *(ii)*.

*(iii):* To prove *(iii)*, let $\mathfrak{q}$ be a non-zero graded prime ideal of $(\mathrm{gr.}A)_S$. Then, by *(ii)*, $\mathfrak{q} \subseteq S^{-1}\mathfrak{p}_i$ for some $i$. But, by *(i)*, $S^{-1}\mathfrak{p}_i$ is a prime ideal of height one in $(\mathrm{gr.}A)_S$, and so we must have $\mathfrak{q} = S^{-1}\mathfrak{p}_i$. This completes the proof of the proposition. □

**Proposition 3.5.** *Every graded ideal in $(\mathrm{gr.}A)_S$ is principal.*

*Proof.* By axiom (C3), the ring $\mathrm{gr.}A$ is factorial and so the localization $(\mathrm{gr.}A)_S$ is also factorial by [B-CA, Chap. VII, §3.4, Prop. 3]. Now take $\mathfrak{b}$ to be any non-zero graded ideal of $(\mathrm{gr.}A)_S$. Thus the module $(\mathrm{gr.}A)_S/\mathfrak{b}$ is graded, and so every element of $\mathrm{Ass}(\mathrm{gr.}A)_S/\mathfrak{b})$ must be a graded prime ideal of $(\mathrm{gr.}A)_S$ [B-CA, Chap. IV, §3.1, Prop. 1]. But Prop. 3.4 shows that every non-zero prime ideal of $(\mathrm{gr.}A)_S$ which is graded has height one. It follows from [B-CA, Chap. VII, §1.6, Prop. 10] therefore that the ideal $\mathfrak{b}$ is divisorial, whence $\mathfrak{b}$ is principal by the definition of a factorial ring. This completes the proof of the proposition. □

While it is technically very convenient that we shall not need to pass to the algebraic microlocalizations of our modules, we shall nevertheless now make essential use of the beginnings of that theory. We recall the definition of the *principal symbol* map. If $x$ is any non-zero element of $M$, we define $\overline{x}$ to be the class of $x$ modulo $F_{n-1}M$, where $n$ is the unique integer such that $x$ belongs to $F_nM$ but not to $F_{n-1}M$. Similarly, we can define by sending $a$ to $\overline{a}$ the map $A \setminus \{0\} \to \text{gr}.A$ . Now let $S$ be the multiplicative set of non-zero homogeneous elements of gr.$A$ defined above in 3.3.

**Definition 3.6.** The set $T$ is defined as $T = \{t \in A : \overline{t} \in S\}$.

As $S$ is a multiplicatively closed set of non-zero elements of gr.$A$, it is plain that $T$ is a multiplicatively closed set of non-zero elements of $A$. What is not at all obvious, and lies at the heart of our proof is the following well-known and basic result. We recall that a multiplicatively closed subset $R$ of non-zero elements of $A$ is a *left and right Ore set* in $A$ if, for each $r$ in $R$ and $a$ in $A$, we have $aR \cap rA$ and $Ra \cap Ar$ are both non-empty sets.

**Proposition 3.7.** *$T$ is a left and right Ore set in $A$.*

*Proof.* We understand that results of this kind go back to Kashiwara [Ka]. For purely algebraic proofs, we refer the reader to [Li], or to the last part of [WK] for a proof due to Björk. □

In view of Prop. 3.7, an entirely parallel construction to that in the commutative case (see [Ste] for instance) proves the existence of the localized domains $T^{-1}A$ and $AT^{-1}$. The left and right localizations are in fact isomorphic, and from now on, we shall denote either of these two by $A_T$. Moreover, as is explained in [Li], $A_T$ is naturally endowed with a separated and exhaustive filtration $F.A_T$. Now, it is also shown in [Li, Prop. 2.3], that this filtration has the property that

$$\text{gr}.(A_T) = (\text{gr}.A)_S.$$

If $N$ is any finitely generated (left) module over $A$ with a good filtration $F.N$, then we define the localization $N_T = A_T \otimes_A N$. This is a finitely generated left module over $A_T$, and comes equipped with a natural good filtration which we denote by $F.N_T$. Further (cf. [Li, Cor. 2.5 and Prop. 2.6]), we have

$$\text{gr}.N_T \simeq (\text{gr}.N)_S.$$

Even though we have imposed the axiom that $A$ is complete with respect to its filtration, it will not in general be true that $A_T$ is complete with respect to the filtration $F.A_T$. However, for our purposes, the following weaker result suffices.

**Proposition 3.8.** *The ring $A_T$ with its natural filtration $F.A_T$ is a Zariski ring.*

*Proof.* See [Li, Prop. 2.8]. □

We can now prove the first main result of this section.

**Proposition 3.9.** *Every left and right ideal of $A_T$ is principal.*

*Proof.* By symmetry, it suffices to prove that any left ideal $L$ of $A_T$ is generated by one element. Endow $L$ with the filtration induced from $A_T$. By a basic property of Zariski rings, (see [LVO, Chap. I, §5, Cor. 5.5]), this induced filtration on $L$ is good. Plainly, gr.$L$ is a graded ideal in gr.$A_T = (\text{gr}.A)_S$. By Prop. 3.5, gr.$L$ is therefore principal, and thus we can find a homogeneous $z$ in gr.$L$ such that gr.$L = (\text{gr}.A)_S \cdot z$. Now pick $w$ to be any element of $L$ such that $\overline{w} = z$. But now, by another basic property of Zariski rings [LVO, Chap. I, §5, Cor. 5.5], we conclude that $L = A_T w$. We stress that, thanks to the remarkable properties of Zariski rings, we are able to carry out this last step without passing to the microlocalization of $A_T$. □

We now consider the left $A_T$-module $M_T$ which is finitely generated and torsion over $A_T$, because $M$ is assumed to be finitely generated and torsion over $A$. Since $A_T$ is a principal ideal domain by Prop. 3.9, the classical theory of finitely generated modules over principal ideal domains (see [Ja, Chap. 3, Thm. 19] or [MCR, Cor. 5.7.19]), which are not necessarily commutative, immediately gives the following result.

**Proposition 3.10.** *There exist elements $w_1, \cdots, w_m$ in $M_T$ such that*

$$M_T = A_T w_1 \oplus \cdots \oplus A_T w_m .$$

Under the present assumptions on the ring $A$ we now are able, using Prop. 3.10, to give the following alternative proof for Cor. 2.5.

**Proposition 3.11.** *Let $M$ be any $A$-torsion module in $\mathfrak{M}(A)$, and let $M_\text{o}$ denote the maximal pseudo-null submodule of $M$. Then there exist non-zero left ideals $L_1, \cdots, L_m$ of $A$ along with an injective $A$-homomorphism*

$$\phi \ : \ \bigoplus_{i=1}^{m} A/L_i \ \longrightarrow \ M/M_\text{o}$$

*such that* $\mathrm{Coker}(\phi)$ *is pseudo-null.*

**Corollary 3.12.** *Let $M$ be as above, and assume, in addition, that $M$ has no non-zero pseudo-null submodule. Then there exist non-zero left ideals $L_1, \cdots, L_m$ of $A$ and an injective $A$-homomorphism*

$$\phi : \bigoplus_{i=1}^{m} A/L_i \longrightarrow M$$

*such that* $\mathrm{Coker}(\phi)$ *is pseudo-null.* □

A crucial step in the proof of these results is given by the following proposition, whose proof will be given at the end of this section. We write $\psi : M \to M_T$ for the canonical $A$-homomorphism, and define $M' := \mathrm{Im}(\psi)$.

**Proposition 3.13.** *Let $N$ be any subquotient module of the torsion module $M$. Then $N_T = 0$ implies that $N$ is pseudo-null. The converse statement holds if $A$ is assumed to be Auslander regular.*

Assuming Prop. 3.13, we now deduce Prop. 3.11 from Prop. 3.10. Put $N := \mathrm{Ker}(\psi)$. Since $\mathrm{Ker}(\psi)$ is precisely the set of $T$-torsion elements in $M$, we have $N_T = 0$ and so $N$ is pseudo-null by Prop. 3.13. Now $M'$ is an $A$-submodule of $M_T$ with $M'_T = M_T$. Let $w_1, \cdots, w_m$ be the elements of $M_T$ given by Prop. 3.10. Since we are clearly free to multiply these elements by any element of $T$, the fact that $M'_T = M_T$ allows us to assume that each of $w_1, \cdots, w_m$ belongs to $M'$. We now define the $A$-submodule $M''$ of $M'$ by

$$M'' := Aw_1 \oplus \cdots \oplus Aw_m ,$$

where the sum is clearly direct because $w_1, \cdots, w_m$ are even linearly independent over $A_T$. But $N := M'/M''$ is a quotient of $M$ with $N_T = 0$. Hence $N$ is pseudo-null by Prop. 3.13. Moreover, the map $a \mapsto aw_i$ induces an isomorphism of $A$-modules from $A/L'_i$ to $Aw_i$, where $L'_i$ is the annihilator of $w_i$ in $A$. The composed map

$$\bigoplus_{i=1}^{m} A/L'_i \xrightarrow{\cong} \bigoplus_{i=1}^{m} Aw_i = M'' \subseteq M' \cong M/\mathrm{Ker}(\psi)$$

is an injective $A$-homomorphism with pseudo-null cokernel. It induces the map $\phi$ in the assertion once we replace each $L'_i$ by the unique left ideal $L_i$ such that $L_i/L'_i$ is the maximal pseudo-null submodule of $A/L'_i$. This completes the proof of Prop. 3.11. □

*Proof of Prop. 3.13.* Take $N$ to be any $A$-subquotient of $M$. Our filtration $F.M$ being fixed, we always endow $N$ with the natural filtration which makes gr.$N$ a subquotient of gr.$M$, namely, if we write $N = V/V'$, where $V' \subset V$ are submodules of $M$, then we take the filtration on $V$ induced from that on $M$, and then the quotient filtration on $N$. By a fundamental property of Zariski rings (see [LVO, Chap. II, §2.1, Thm. 2]), this filtration on $N$ is good because we start with a good filtration on $M$. We next observe that since gr.$N$ is a subquotient of gr.$M$, the support of gr.$N$ must be a subset of the support of gr.$M$, whence (see Remark 3.2)

$$W(N) \subseteq W(M).$$

Finally, it will be important for us to distinguish between homogeneous localization and ordinary localization of a gr.$A$-module at a prime ideal of gr.$A$. Let $R$ denote a gr.$A$-module, and $\mathfrak{p}$ a prime ideal of gr.$A$. As usual, we write $R_\mathfrak{p}$ for the localization of $R$ with respect to the multiplicative set $S_\mathfrak{p} = \text{gr}.A \setminus \mathfrak{p}$. If in addition, $R$ is a graded gr.$A$-module and $\mathfrak{p}$ is a graded prime ideal, we shall write $R_{\mathfrak{p}^*}$ for the localization of $R$ with respect to the multiplicative set $S_{\mathfrak{p}^*}$ of all homogeneous elements of gr.$A \setminus \mathfrak{p}$.

Suppose now that $N$ is a subquotient of $M$ with $N_T = 0$. We must show that $N$ is pseudo-null. Since $A$ is Zariski, it is well-known that (see [LVO, Chap. III, §2.2])

$$j_A(N) \geq j_{\text{gr}.A}(\text{gr}.N),$$

where $j$ denotes the grade of the relevant module as defined in §1. Thus, by Lemma 1.1, to prove that $N$ is pseudo-null, it suffices to prove that gr.$N$ is pseudo-null as a gr.$A$-module. Hence, by Example 1.3, we must prove that $(\text{gr}.N)_\mathfrak{p} = 0$ for every prime ideal $\mathfrak{p}$ in $W(N)$. But, as remarked in the previous paragraph, $W(N) \subseteq W(M)$ because $N$ is a subquotient of $M$, and so it suffices to show that $(\text{gr}.N)_{\mathfrak{p}_i} = 0$ for $1 \leq i \leq r$. Now our good filtration on $N$ gives rise (see [Li, Cor. 2.5]) to a good filtration on the $A_T$-module $N_T$, with the property that

$$\text{gr}.N_T = (\text{gr}.N)_S.$$

Hence we must have $(\text{gr}.N)_S = 0$. But this automatically implies that for $1 \leq i \leq r$, we have $(\text{gr}.N)_\mathfrak{p} = 0$, because $S \subseteq S_{\mathfrak{p}_i^*} \subseteq S_{\mathfrak{p}_i}$.

Conversely, assume that $N$ is pseudo-null, and also that $A$ is Auslander regular. This last assumption guarantees that we have $j_A(N) = j_{\text{gr}.A}(\text{gr}.N)$ (see [LVO, Chap. III, §2.2, Thm. 5]). Hence, by Lemma 1.4, gr.$N$ is a pseudo-null gr.$A$-module, since $N$ is pseudo-null. Let $\mathfrak{p}$ denote any graded prime ideal of height 1 of gr.$A$. As gr.$N$ is pseudo-null, we have $(\text{gr}.N)_\mathfrak{p} = 0$. In fact, we claim that we have the stronger assertion that $(\text{gr}.N)_{\mathfrak{p}^*} = 0$. Let $w_1, \cdots, w_k$ denote a set of homogeneous generators for gr.$N$ as a gr.$A$-module. For each $j = 1, \cdots, k$, we can find $s_j \in S_\mathfrak{p}$ such that $s_j w_j = 0$. Writing $s_j$ as a sum of its homogeneous components $s_{j,h}$, we must always have $s_{j,h} w_j = 0$. But at least one of the $s_{j,h}$ does not belong to $\mathfrak{p}$. Hence we must have $s_{j,h} \in S_{\mathfrak{p}^*}$ for some $h$, which proves that $(\mathrm{gr}.N)_{\mathfrak{p}^*} = 0$. We claim that this, in turn, implies that $(\mathrm{gr}.N)_S = 0$. Indeed, the prime ideals in $W(M) = \{\mathfrak{p}_1, \cdots, \mathfrak{p}_r\}$ are all graded, and the previous argument has shown that there exists a homogeneous element $t_i$ in $\mathrm{ann}_{\mathrm{gr}.A}(\mathrm{gr}.N)$ such that $t_i$ does not belong to $\mathfrak{p}_i$. Let $\mathfrak{a}$ be the graded ideal of $\mathrm{gr}.A$ generated by $t_1, \cdots, t_r$. By the usual argument, the ideal $\mathfrak{a}$ cannot be contained in the union $\mathfrak{p}_1 \cup \cdots \cup \mathfrak{p}_r$, as this would imply that $\mathfrak{a} \subseteq \mathfrak{p}_i$ for some $i$. Hence, picking a homogeneous element $t$ of $\mathfrak{a}$ which does not belong to $\mathfrak{p}_1 \cup \cdots \cup \mathfrak{p}_r$, we have $t \cdot \mathrm{gr}.N = 0$ and $t \in S$, so that $(\mathrm{gr}.N)_S = 0$, as asserted. Recall that our goal is to show that $N_T = 0$. Now, as always, our good filtration on $N$ gives rise to a good filtration on the $A_T$-module $N_T$ with the property that $\mathrm{gr}.N_T = (\mathrm{gr}.N)_S$. Hence we have $\mathrm{gr}.N_T = 0$. But $A_T$ is a Zariski ring, and thus every good filtration is separated [LVO, Chap. II, §2.1, Thm.2], and so $\mathrm{gr}.N_T = 0$ implies that $N_T = 0$. This completes the proof of Prop. 3.13. □

**Remark 3.14.** We remark that a Zariski ring $A$ with the property that $\mathrm{gr}.A$ is Auslander regular is itself Auslander regular [LVO, Chap. III, §2.2]. Then the above proposition shows that for the dimension filtration on $M$ (cf. §1) the submodule $\Delta^2(M)$ is precisely the $T$-torsion submodule of $M$.

Finally, we record the following corollary that emerges from our methods, and which is of independent interest.

**Corollary 3.15.** *Let $A$ be a Zariski ring such that $\mathrm{gr}.A$ is a polynomial ring in finitely many variables over a field $k$. Suppose $M \in \mathfrak{M}(A)$ is a torsion module and $T$ is defined as before. Assume further that $M$ has no $T$-torsion. Then there exists a good filtration on $M$ such that $\mathrm{gr}.M$ is $S$-torsion free.*

*Proof.* Indeed, $A$ is then Auslander regular, and if $M$ has no $T$-torsion, then by Prop. 3.13 $M$ is pure [LVO, Chap. III, §4.2, Thm. 6]. It is well-known that in this case there is a good filtration [LVO, Chap. III, §4.2, Thm. 13 (4)] on $M$ such that $\mathrm{gr}.M$ is a pure module over $\mathrm{gr}.A$. But this means in particular [LVO, Chap. III, §4.3, Prop. 5] that $\mathrm{Ass}(\mathrm{gr}.M) = W(M)$. Recall that the set of zero divisors on $\mathrm{gr}.M$ is equal to the union of the prime ideals in $\mathrm{Ass}(\mathrm{gr}.M)$. By definition, $S$ has trivial intersection therefore with the set of zero divisors and hence $\mathrm{gr}.M$ has no $S$-torsion. □

## 4. The characteristic ideal

Let us return to the general situation when $A$ is a noetherian maximal order without zero divisors and with skew-field of fractions $D$. We always assume that $A \neq D$. There are certain additional features of the structural result Prop. 2.3 which at least partly generalize [B-CA, Chap. VII §4.5] to the non-commutative situation. As before the annihilator ideal of an $A$-module $M$ is denoted by $\mathrm{ann}_A(M)$. If $\mathrm{ann}_A(M) \neq 0$ then the module $M$ is called *bounded*. The annihilator ideal of an object $\mathcal{M}$ in the quotient category $\mathrm{Mod}(A)/\mathcal{C}_A^1$ is defined by

$$\mathrm{ann}(\mathcal{M}) := \sum \{\mathrm{ann}_A(M) : Q(M) \cong \mathcal{M}\} \ .$$

The object $\mathcal{M}$ is called *completely faithful* if $\mathrm{ann}(\mathcal{N}) = 0$ for any non-zero subquotient object $\mathcal{N}$ of $\mathcal{M}$. It is called *(locally) bounded* if $\mathrm{ann}(\mathcal{N}) \neq 0$ for any subobject $\mathcal{N} \subseteq \mathcal{M}$ (of finite type). Since $A$ is noetherian it follows from [Ste, Chap. XIII, Prop. 1.1] that $\mathcal{N}$ is of finite type if and only if it is of the form $\mathcal{N} \cong Q(N)$ for some finitely generated module $N$.

Under the condition that $\mathrm{K}_1\mathrm{dim}(A) \leq 1$ there is an alternative description of $\mathrm{ann}(Q(M))$ for any finitely generated $A$-torsion module $M$ in terms of the left ideals $L_1, \ldots, L_m$ occuring in Cor. 2.5. For $i = 1, \ldots, m$, let $J_i$ be the maximal two sided ideal of $A$ contained in $L_i$, and let $J := \bigcap_{i=1}^m J_i$. Then we have

$$\mathrm{ann}(Q(M)) = J \ .$$

This follows easily on combining Cor. 2.5 with the observation that $\mathrm{ann}_A(A/L_i) = J_i$, and using the fact [Rob, Lemma 2.5] that $\mathrm{ann}(Q(M)) = \mathrm{ann}_A(M/M_\mathrm{o})$, where again $M_\mathrm{o}$ denotes the maximal pseudo-null submodule of $M$. Note also that $M$ is bounded, or equivalently $J \neq 0$ if and only if all of $J_1, \ldots, J_m$ are non-zero.

**Proposition 4.1.** *(i) Any object $\mathcal{M}$ in $\mathcal{C}_A^0/\mathcal{C}_A^1$ decomposes uniquely into a direct sum $\mathcal{M} = \mathcal{M}_0 \oplus \mathcal{M}_1$ where $\mathcal{M}_0$ is completely faithful and $\mathcal{M}_1$ is locally bounded. (ii) Any completely faithful object of finite length in $\mathrm{Mod}(A)/\mathcal{C}_A^1$ is cyclic.*

*Proof.* See [Ch2, Prop. 4.2.2 and Cor. 4.2.3]. □

Let $(\mathcal{C}_A^0/\mathcal{C}_A^1)^b$ denote the full subcategory of all bounded objects in $\mathcal{C}_A^0/\mathcal{C}_A^1$. As a consequence of [Rob, Lemma 2.4] this category is an abelian subcategory. We want to investigate the Grothendieck group $G((\mathcal{C}_A^0/\mathcal{C}_A^1)^b)$ of objects of finite length in the category $(\mathcal{C}_A^0/\mathcal{C}_A^1)^b$ modulo short exact sequences. By the Jordan-Hölder theorem $G((\mathcal{C}_A^0/\mathcal{C}_A^1)^b)$ is the free abelian group on the set of isomorphism classes of simple objects in the category $(\mathcal{C}_A^0/\mathcal{C}_A^1)^b$. The tool we will use is a certain group $G(A)$ of fractional ideals of the maximal order $A$ whose construction we now briefly review, after first introducing more notation.

For any left or right $A$-module $M$ let $M^* := \operatorname{Hom}_A(M, A)$ denote the dual right or left $A$-module.

A non-zero left (right) $A$-submodule $L \subseteq D$ is called a fractional left (right) ideal if there is a $v \in D^\times$ such that $L \subseteq Av$ ($L \subseteq vA$). If $I \subseteq D$ is a fractional left as well as right ideal we call it a fractional ideal. For any fractional left ideal $L$ the map
$$\{x \in D : Lx \subseteq A\} \xrightarrow{\cong} L^*$$
$$x \longmapsto [b \mapsto bx]$$
is an isomorphism of right $A$-modules. In a maximal order one moreover has
$$\{x \in D : Lx \subseteq A\} = \{x \in D : LxL \subseteq L\} =: L^{-1} .$$

Defining
$$\overline{L} := (L^{-1})^{-1}$$
we have an isomorphism $\overline{L} \cong L^{**}$. It follows from the left-right symmetry of the definition of $L^{-1}$ that a fractional ideal $I \subseteq D$ is reflexive as a left $A$-module if and only if it is reflexive as a right $A$-module. A reflexive fractional ideal, resp. a non-zero reflexive two sided ideal in $A$, is called a fractional c-ideal, resp. a c-ideal. According to Asano ([Asa], compare also [Ch1],[Mar], or [AKMU, Lemma 1.1]) one has:

– If $I$ is a fractional ideal then $\overline{I}$ is a fractional c-ideal.

– The set of all fractional c-ideals is an abelian group $G(A)$ with respect to the product
$$I \cdot J := \overline{IJ}$$
and the inverse $I \longmapsto I^{-1}$.

– Every maximal element in the set of all c-ideals is a prime ideal of height one (i.e., a minimal non-zero prime ideal). The group $G(A)$ is the free abelian group on the set $\mathcal{P} = \mathcal{P}(A)$ of all prime c-ideals.

– For any pairwise distinct $P_1, \ldots, P_r \in \mathcal{P}$ and any $n_1, \ldots, n_r \in \mathbb{N}$ one has
$$P_1^{n_1} \cdot \ldots \cdot P_r^{n_r} = P_1^{n_1} \cap \ldots \cap P_r^{n_r} .$$

If $C(P)$, for any $P \in \mathcal{P}$, denotes the multiplicative set of elements in $A$ which are regular modulo $P$ then it is proved in [Ch0, Prop. 1.8(b), 1.9, and 2.5] or [Ch1, Prop. 3.3.4 and Lemma 3.3.5] that:

– $C(P)$ satisfies the left and right Ore condition.

– The localization $A_{C(P)} = {}_{C(P)}A$ is a bounded left and right principal maximal order with Jacobson radical $\mathrm{Jac}(A_{C(P)}) = PA_{C(P)}$ such that $A_{C(P)}/PA_{C(P)}$ is simple artinian and such that the non-zero two sided ideals in $A_{C(P)}$ are precisely the powers of $\mathrm{Jac}(A_{C(P)})$.

– For any non-zero $a \in A$ there are at most finitely many $P \in \mathcal{P}$ such that $a \notin C(P)$.

**Lemma 4.2.** *For any $A$-module $M$ in $\mathcal{C}_A^0$ the following assertions are equivalent:*
*(i) $Q(M)$ is completely faithful*
*(ii) $A_{C(P)} \otimes_A M = 0$ for any $P \in \mathcal{P}$.*

*Proof.* See [Ch2, Lemma 4.2.1]. □

In particular, this lemma implies that $A_{C(P)} \otimes_A N = 0$ for every pseudo-null $A$-module $N$. Since localization is exact we therefore obtain, for each $P \in \mathcal{P}$, a well-defined exact functor

$$\mathcal{C}_A^0/\mathcal{C}_A^1 \longrightarrow \mathrm{Mod}(A_{C(P)})$$

$$\mathcal{M} = Q(M) \longmapsto A_{C(P)} \otimes_A \mathcal{M} := A_{C(P)} \otimes_A M.$$

Suppose now that $\mathcal{M}$ is an object of finite length in $(\mathcal{C}_A^0/\mathcal{C}_A^1)^b$. As remarked earlier we then have $\mathcal{M} \cong Q(M)$ for some finitely generated torsion $A$-module $M$. Hence $A_{C(P)} \otimes_A \mathcal{M} \cong A_{C(P)} \otimes_A M$ is a finitely generated torsion $A_{C(P)}$-module. Since $A_{C(P)}$ is principal this latter module is of finite length $\ell_P(\mathcal{M})$. The last fact which we recalled before Lemma 4.2 implies that $\ell_P(\mathcal{M}) \neq 0$ for at most finitely many $P$. Hence

$$\chi(\mathcal{M}) := \prod_{P \in \mathcal{P}} P^{\ell_P(\mathcal{M})}$$

is a well-defined c-ideal. We call it the *characteristic ideal* of the object $\mathcal{M}$. It follows from Lemma 4.2 that $\chi(\mathcal{M}) = A$ implies $\mathcal{M} = 0$.

**Lemma 4.3.** *Let $M$ be a finitely generated bounded torsion $A$-module, and let $M_o$ be the maximal pseudo-null submodule of $M$. We then have:*
*(i) $\mathrm{ann}_A(M/M_o) = \mathrm{ann}(Q(M))$ is a c-ideal,*
*(ii) $\chi(Q(M)) \subseteq \mathrm{ann}(Q(M))$,*
*(iii) $\chi(Q(M))$ and $\mathrm{ann}(Q(M))$ have the same prime factors.*

*Proof.* Replacing $M$ by $M/M_{\mathrm{o}}$, we may assume that $M$ has no non-zero pseudo-null submodules.

*(i):* According to [Ch2, Cor. of Thm. 2.5] and [Rob, Lemma 2.5] we have $I := \mathrm{ann}_A(M) = \mathrm{ann}(Q(M))$. Hence $I \neq 0$, by assumption, and it remains to show that $\overline{I} \subseteq I$. But

$$I \cdot I^{-1} \cdot \overline{I} = I \cdot I^{-1} \cdot (I^{-1})^{-1} \subseteq \mathrm{ann}_A(M) \ .$$

Since the ideal $J := I \cdot I^{-1}$, by [Ch2, Lemma 2.2], is such that the quotient $A/J$ is pseudo-null the submodule $\overline{I} \cdot M$ being annihilated by $J$ is pseudo-null as well. It follows that $\overline{I} \cdot M = 0$ which means that $\overline{I} \subseteq \mathrm{ann}_A(M) = I$.

*(ii):* Any $A_{C(P)}$-module of length $m$ is annihilated by $\mathrm{Jac}(A_{C(P)})^m$. Hence for any $P \in \mathcal{P}$ we have

$$A_{C(P)} \otimes_A \chi(Q(M))M \subseteq A_{C(P)} \otimes_A P^{\ell_P(Q(M))}M$$

$$\subseteq \mathrm{Jac}(A_{C(P)})^{\ell_P(Q(M))}(A_{C(P)} \otimes_A M) = 0 \ .$$

Lemma 4.2 then implies that $\chi(Q(M))M$ is a pseudo-null submodule of $M$ and hence is zero by our assumption on $M$.

*(iii):* The prime factors of $\mathrm{ann}_A(M)$ are precisely the prime c-ideals which contain $\mathrm{ann}_A(M)$. Because of *(ii)* it therefore suffices to show that any $P \in \mathcal{P}$ such that $A_{C(P)} \otimes_A M \neq 0$ contains $\mathrm{ann}_A(M)$. Since $\mathrm{ann}_A(M)$ is reflexive by *(i)* it follows from [Ch2, Lemma 3.3] that $\mathrm{ann}_A(M) \subseteq P$ if and only if $\mathrm{ann}_A(M) \cap C(P) = \emptyset$. The latter condition certainly is satisfied if $A_{C(P)} \otimes_A M \neq 0$. □

Thanks to the exactness of localization and the additivity of the length function the formation of the characteristic ideal induces a homomorphism of abelian groups

$$\chi : G((\mathcal{C}_A^0/\mathcal{C}_A^1)^b) \longrightarrow G(A) \ .$$

To understand this homomorphism consider any simple object $\mathcal{M}$ in $(\mathcal{C}_A^0/\mathcal{C}_A^1)^b$. By Prop. 2.3 it is cyclic and therefore a quotient object of $Q(A)$. Since the subobjects of $Q(A)$ are in bijection with the $\mathcal{C}_A^1$-closed left ideals we see that $\mathcal{M}$, up to isomorphism, is of the form

$$\mathcal{M} = Q(A/L)$$

where the left ideal $L$ is maximal among all $\mathcal{C}_A^1$-closed left ideals. In addition, since $\mathrm{ann}_A(A/L) = \mathrm{ann}(\mathcal{M})$, by [Ch2, Cor. of Thm. 2.5] and [Rob, Lemma 2.5], and the latter is non-zero by assumption, the left ideal $L$ is bounded.

Applying [Ch0, Lemma 2.7] we obtain that $L$ is reflexive. By Lemma 4.2 there exists a $P \in \mathcal{P}$ such that $A_{C(P)} \otimes_A \mathcal{M} = A_{C(P)} \otimes_A A/L \neq 0$ and hence that $\mathrm{ann}_A(A/L) \cap C(P) = \emptyset$. Since both $\mathrm{ann}_A(A/L)$ (see Lemma 4.3(i)) and $L$ are reflexive and bounded, we can apply [Ch2, Lemma 3.3] twice (observing that $\mathrm{ann}_A(A/L)$ is the largest two sided ideal contained in $L$) and conclude that $\mathrm{ann}_A(A/L) \subseteq P$ and $L \cap C(P) = \emptyset$. Suppose that $L \subseteq L' \subseteq A$ is another left ideal such that $L' \cap C(P) = \emptyset$. Then $A_{C(P)} \otimes_A A/L' = A_{C(P)}/A_{C(P)}L'$ is non-zero. Hence, by Lemma 4.2 again, $Q(A/L')$ is a non-zero quotient of $\mathcal{M}$ and therefore must be equal to $\mathcal{M}$. This implies that $L'/L$ is pseudo-null. But since $L$ is $\mathcal{C}_A^1$-closed the quotient $A/L$ contains no non-zero pseudo-null submodules. This shows that $L' = L$ and establishes that $L$ is maximal with respect to the property $L \cap C(P) = \emptyset$. Such left ideals are called $P$-*critical* and are studied in [LM]. In particular, it is shown in [LM, Thm. 3.5] that any $P$-critical left ideal contains $P$. In our case this means that

$$\mathrm{ann}(\mathcal{M}) = P \ .$$

Moreover, as explained in the proof of [Ch1, Prop. 3.4.6], if $L$ is $P$-critical then $L_{C(P)} = A_{C(P)}L$ is a maximal left ideal in $A_{C(P)}$ and $L_{C(Q)} = A_{C(Q)}L = A_{C(Q)}$ for any $Q \in \mathcal{P}$ different from $P$. On the one hand this implies that

$$\chi(\mathcal{M}) = P \ .$$

On the other hand it also follows that if we give ourselves a $P' \in \mathcal{P}$ and choose a $P'$-critical left ideal $L'$ then $\chi(Q(A/L')) = P'$. By Lemma 4.2, the object $Q(A/L')$ in $(\mathcal{C}_A^0/\mathcal{C}_A^1)^b$ necessarily is simple.

Suppose now that $\mathcal{N}$ is a second simple object in $(\mathcal{C}_A^0/\mathcal{C}_A^1)^b$ such that $\mathrm{ann}(\mathcal{N}) = \mathrm{ann}(\mathcal{M}) = P$. Let $\mathcal{N}$ be of the form $\mathcal{N} = Q(A/L')$ with $L'$ maximal among the $\mathcal{C}_A^1$-closed left ideals. Then $L'$ also is $P$-critical and [LM, Lemma 3.7] says that there are elements $a \in A \setminus L$ and $b \in A \setminus L'$ such that

$$a^{-1}L := \{x \in A : xa \in L\} = \{y \in A : yb \in L'\} =: b^{-1}L' \ .$$

Consider the injective map

$$A/a^{-1}L \xrightarrow{\cdot a} A/L \ .$$

Since $A/L$ has no non-zero pseudo-null submodules $Q(A/a^{-1}L)$ is a non-zero subobject of $\mathcal{M}$ and therefore is isomorphic to $\mathcal{M}$. We obtain

$$\mathcal{M} \cong Q(A/a^{-1}L) = Q(A/b^{-1}L') \cong \mathcal{N} \ .$$

Altogether this proves the following result.

**Proposition 4.4.** *The homomorphism* $\chi : G((\mathcal{C}_A^0/\mathcal{C}_A^1)^b) \xrightarrow{\cong} G(A)$ *is an isomorphism. More precisely, $\chi$ induces a bijection between the set of isomorphism classes of simple objects in $(\mathcal{C}_A^0/\mathcal{C}_A^1)^b$ and the set $\mathcal{P}(A)$ of prime c-ideals.* □

Suppose that $A$ is a maximal order satisfying $K_1\dim(A) \leq 1$. We then can reformulate the above result as follows. Let $\mathcal{C}_A$ denote the abelian category of all finitely generated torsion $A$-modules with non-zero annihilator ideal and let $\mathcal{Z}_A \subseteq \mathcal{C}_A$ denote the subcategory of pseudo-null modules. The isomorphism classes of simple objects in the quotient category $\mathcal{C}_A/\mathcal{Z}_A$ are, via the characteristic ideal, in bijection with the set $\mathcal{P}$.

Following [AKMU] we call $A$ a *unique factorization ring* if any prime c-ideal in $A$ is principal. Here a two sided ideal $I \subseteq A$ is called principal if $I = Ab$ (and then necessarily also $I = bA$ by [Ch1, Lemma 2.2.9]) for some $b \in A$. It is clear then that all fractional c-ideals of $A$ are principal. In this case the characteristic ideal $\chi(\mathcal{M})$ determines, through its generator, a *characteristic element* for the object $\mathcal{M}$ in $A$ which is well-defined up to a unit in $A$.

## 5. Finitely generated modules

Under the same assumptions as in the previous section we now want to look at the structure of an arbitrary finitely generated $A$-module $M$ modulo pseudo-null modules. Our aim is to at least partly generalize [B-CA, Chap. VII, §4.4, Example (3) and Thm. 4] which says that, if $A$ is commutative, then the torsion submodule of $M$ always splits off in the quotient category $\mathrm{Mod}(A)/\mathcal{C}_A^1$.

**Lemma 5.1.** *Let $P \in \mathcal{P}(A)$ and let $M$ be an $A$-module such that $Q(M) \cong Q(A/L)$ for some $P$-critical left ideal $L \subseteq A$. The natural map $M \to A_{C(P)} \otimes_A M$ viewed as a homomorphism of $A$-modules induces an isomorphism $Q(M) \xrightarrow{\cong} Q(A_{C(P)} \otimes_A M)$.*

*Proof.* We may assume that $M = A/L$. From the proof of Prop. 4.4 we know the following facts:
(a)   $L$ contains $P$,
(b)   $A_{C(Q)} \otimes_A M = A_{C(Q)}/A_{C(Q)}L = 0$ for any $Q \in \mathcal{P}$ different from $P$.
Consider now the natural map
$$M = A/L \longrightarrow A_{C(P)} \otimes_A M = A_{C(P)}/A_{C(P)}L \ .$$
Because of (a) the $A$-modules on both sides are bounded. For our claim that this map induces an isomorphism in the quotient category it therefore suffices, by Lemma 4.2, to show that the localized maps
$$A_{C(Q)} \otimes_A M \longrightarrow A_{C(Q)} \otimes_A A_{C(P)} \otimes_A M$$
are isomorphisms for any $Q \in \mathcal{P}$. For $Q = P$ this is obvious. Suppose therefore that $Q \neq P$. Then the left hand side is zero by (b). Since $A_{C(Q)} \otimes_A A_{C(P)} = D$ by [Ch1, Lemma 3.5.2] the right hand side is zero as well, because $M$ is $A$-torsion. $\square$

**Proposition 5.2.** *Let $N$ be a finitely generated $A$-torsion free $A$-module, and let $\mathcal{M}$ be a bounded object of finite length in the quotient category $\mathfrak{A} := \mathrm{Mod}(A)/\mathcal{C}_A^1$; then*
$$\mathrm{Ext}^1_{\mathfrak{A}}(Q(N), \mathcal{M}) = 0 \ .$$

*Proof.* First of all we note that $\mathcal{M}$ being bounded lies in the subcategory $\mathcal{C}_A^0/\mathcal{C}_A^1$. By induction with respect to the length of $\mathcal{M}$ we may assume that $\mathcal{M}$ is simple. In the proof of Prop. 4.4 we have seen that there is then a prime c-ideal $P \in \mathcal{P}$ such that $\mathcal{M} \cong Q(A/L)$ for some $P$-critical left ideal $L \subseteq A$. We therefore have to show that any given extension

(1) $$0 \longrightarrow Q(A/L) \longrightarrow \mathcal{E} \longrightarrow Q(N) \longrightarrow 0$$

in $\mathfrak{A}$ splits. We choose an exact sequence of $A$-modules

(2) $$0 \longrightarrow M \longrightarrow E \longrightarrow N' \longrightarrow 0$$

such that (1) is isomorphic to the exact sequence arising from (2) by applying the functor $Q$ (see [Rob, Lemma 1.1]). Consider now the exact sequence

($2_P$) $$0 \longrightarrow A_{C(P)} \otimes_A M \longrightarrow A_{C(P)} \otimes_A E \longrightarrow A_{C(P)} \otimes_A N' \longrightarrow 0 \ .$$

The $A_{C(P)}$-module $A_{C(P)} \otimes_A M \cong A_{C(P)} \otimes_A A/L$ is finitely generated and torsion whereas the $A_{C(P)}$-module $A_{C(P)} \otimes_A N' \cong A_{C(P)} \otimes_A N$ is finitely generated and torsion free. Since $A_{C(P)}$ is a (non-commutative) principal ideal domain it follows [MCR, Lemma 5.7.4] that the sequence ($2_P$) splits. Let $\sigma : A_{C(P)} \otimes_A E \longrightarrow A_{C(P)} \otimes_A M$ be a splitting of ($2_P$). If $\sigma_{\mathrm{o}}$ denotes the restriction of $\sigma$ to $E$ then we deduce from Lemma 5.1 that $Q(\sigma_{\mathrm{o}})$ induces a splitting of (1). $\square$

**Remark 5.3.** *Any bounded object $\mathcal{M}$ in the quotient category $\mathrm{Mod}(A)/\mathcal{C}_A^1$ has injective dimension $\leq 1$.*

*Proof.* Let $\mathcal{M} = Q(M)$ and let

$$0 \longrightarrow M \longrightarrow E_0 \longrightarrow E_1 \longrightarrow \ldots$$

be the "minimal" injective resolution of the $A$-module $M$ as in §1. With $M$ all the injective $A$-modules $E_i$ are $A$-torsion [Ste, Chap. VI, Cor. 6.8, Prop. 7.3]. Hence the exact sequence

$$0 \longrightarrow \mathcal{M} \longrightarrow Q(E_0) \longrightarrow Q(E_1) \longrightarrow \ldots$$

lies in $\mathcal{C}_A^0/\mathcal{C}_A^1$. According to Prop. 4.1*(i)* each

$$Q(E_i) = Q(E_i)_0 \oplus Q(E_i)_1$$

decomposes uniquely into a completely faithful object $Q(E_i)_0$ and a locally bounded object $Q(E_i)_1$. Since $\mathcal{M}$ is bounded the sequence

$$0 \longrightarrow \mathcal{M} \longrightarrow Q(E_0)_1 \longrightarrow Q(E_1)_1 \longrightarrow \ldots$$

still is exact. We claim that $Q(E_i)_1$ is injective in $\mathrm{Mod}(A)/\mathcal{C}_A^1$ for any $i \geq 0$ and is equal to zero for $i \geq 2$. By [Gab, Chap. III, §3, Cor. 2] we have, for any $i \geq 0$, another decomposition

$$Q(E_i) = Q(E_i') \oplus \mathcal{E}_i$$

where $E_i'$ is the injective hull of some pseudo-null $A$-module and $\mathcal{E}_i$ is an injective object in $\mathrm{Mod}(A)/\mathcal{C}_A^1$. At this point we have to introduce the intermediate ring $A \subseteq A_0 \subseteq D$ defined by

$$A_0 := \{u \in D : uI \subseteq A \text{ for some non-zero two sided ideal } I \subseteq A\} \ .$$

Let $\mathcal{C}_A^{1/2}$ denote the hereditary torsion theory cogenerated by the injective $A$-module $D \oplus E(A_0/A)$ (cf. §1). Since $E(A_0/A)$ embeds into $E(D/A)$ we have

$$\mathcal{C}_A^1 \subseteq \mathcal{C}_A^{1/2} \subseteq \mathcal{C}_A^0 \ .$$

According to [BHL, Remark 4.4(i)] an $A$-module $N$ lies in $\mathcal{C}_A^{1/2}$ if and only if $A_{C(P)} \otimes_A N = 0$ for any $P \in \mathcal{P}(A)$. It therefore follows from Lemma 4.2 that an $A$-torsion module $N$ lies in $\mathcal{C}_A^{1/2}$ if and only if $Q(N)$ is completely faithful. The subcategory $\mathcal{C}_A^{1/2}$ of $\mathrm{Mod}(A)$ has the important additional property [BHL, Thm. 4.2] that it is stable with respect to the formation of injective hulls. In our situation this shows that

$$Q(E_i') \subseteq Q(E_i)_0$$

and hence that $Q(E_i)_1$ is injective for any $i \geq 0$. Moreover, by [BHL, Remark 4.4(ii)], the sequence

$$0 \longrightarrow A_{C(P)} \otimes_A M \longrightarrow A_{C(P)} \otimes_A E_0 \longrightarrow A_{C(P)} \otimes_A E_1 \longrightarrow \ldots$$

is, for any $P \in \mathcal{P}$, a "minimal" injective resolution of the $A_{C(P)}$-module $A_{C(P)} \otimes_A M$. Since $A_{C(P)}$ as a principal ideal domain has global injective dimension $\leq 1$ we obtain that $A_{C(P)} \otimes_A E_i = 0$ for any $i \geq 2$ and any $P \in \mathcal{P}$. This means, by Lemma 4.2, that $Q(E_i)$ is completely faithful and hence that $Q(E_i)_1 = 0$ for any $i \geq 2$, and completes the proof of our claim. We mention that in the quotient category $\mathrm{Mod}(A)/\mathcal{C}_A^{1/2}$ every object has injective dimension $\leq 1$ [BHL, Remark 4.4(iii)]. $\square$

**Proposition 5.4.** *Let $M$ be a finitely generated $A$-module, and denote by $M_t \subseteq M$ the maximal $A$-torsion submodule. Suppose that $\mathrm{K}_1\dim(A) \leq 1$ and that $Q(M_t)$ does not contain any non-zero completely faithful subobject. Then*

$$Q(M) \cong Q(M_t) \bigoplus Q(M/M_t) \ .$$

*If, in addition, $A$ is Auslander regular then $Q(M/M_t) \cong Q(N)$ for some finitely generated reflexive $A$-module $N$.*

*Proof.* By our assumptions and Prop. 4.1(i) the object $Q(M_t)$ in $\mathrm{Mod}(A)/\mathcal{C}_A^1$ is bounded and of finite length and the $A$-module $M/M_t$ is finitely generated and $A$-torsion free. It therefore follows from Prop. 5.2 that

$$Q(M) \cong Q(M_t) \bigoplus Q(M/M_t) \ .$$

The canonical map $M \longrightarrow M^{**}$ has kernel $M_t$ [MCH, Prop. 3.4.7], and the $A$-module $M^{**}$ is finitely generated and reflexive. If $A$ is Auslander regular then the cokernel of this map is pseudo-null [Ven0, Cor. 1.5.7]. □

## 6. Iwasawa algebras

Let $p$ be a prime number, and let $G$ be a compact $p$-adic Lie group. We recall that the Iwasawa algebra $\Lambda(G)$ is defined to be the completed group ring

$$\Lambda(G) := \varprojlim_{N} \mathbb{Z}_p[G/N]$$

where $N$ runs over the open normal subgroups of $G$. It seems fair to say that, until very recently, Iwasawa algebras have been largely neglected as concrete examples of non-commutative rings, even though they occur very naturally in many questions in arithmetic geometry. Under the assumption that $G$ is $p$-valued in the sense of [Laz], we now show that $\Lambda(G)$ carries an exhaustive and complete filtration, and satisfies the axioms imposed on the ring $A$ in the earlier parts of the paper.

**Definition 6.1.** ([Laz, Chap. III, Def. 2.1.2]) A *$p$-valuation* on $G$ is a function $\omega : G \to (0, \infty]$ satisfying the following axioms for all $g$ and $h$ in $G$:
1) $\omega(1) = \infty$, and $1/(p-1) < \omega(g) < \infty$ for $g \neq 1$;
2) $\omega(gh^{-1}) \geq \min\{\omega(g), \omega(h)\}$;
3) $\omega(g^{-1}h^{-1}gh) \geq \omega(g) + \omega(h)$;
4) $\omega(g^p) = \omega(g) + 1$.

We say that $G$ is *p-valued* if it possesses a $p$-valuation. As B. Totaro emphasized to us, $p$-valued compact $p$-adic Lie groups $G$ possess many useful properties. If $G$ is $p$-valued, it is automatically pro-$p$ and has no element of order $p$. Any closed subgroup of a $p$-valued group $G$ is also $p$-valued. The classic example of a $p$-valued group is the subgroup of $GL_n(\mathbb{Z}_p)$ consisting of all matrices which are congruent to the identity modulo $p$ (resp. modulo 4) if $p$ is odd (resp. if $p = 2$). More generally, if $p > n + 1$, Lazard (see [Laz, p. 101] and [Se3, Chap. I §1.4]) has proven that every closed pro-$p$ subgroup of $GL_n(\mathbb{Z}_p)$ is $p$-valued. The following result, whose proof we shall give below, is essentially contained in Lazard [Laz], and we are grateful to B. Totaro for pointing it out to us.

**Proposition 6.2.** *Assume that $G$ is a compact $p$-adic Lie group, which is $p$-valued. Then $\Lambda(G)$ possesses a complete, separated and exhaustive filtration $F.\Lambda(G)$ such that $\mathrm{gr}.\Lambda(G)$ is isomorphic to the polynomial ring $\mathbb{F}_p[X_0, \cdots, X_d]$ in $d+1$ variables, where $d = \dim(G)$. In particular $\Lambda(G)$ is a noetherian Auslander regular maximal order without zero divisors, such that $\mathrm{K}_1\dim(\Lambda(G)) \leq 1$, and satisfies axioms (C1), (C2) and (C3) of section 3.*

**Remark 6.3.** In the isomorphism from $\mathrm{gr}.\Lambda(G)$ to $\mathbb{F}_p[X_0, \cdots, X_d]$ induced by the filtration $F.\Lambda(G)$, the degree of $X_i$ ($0 \leq i \leq d$) is not necessarily the naive one, namely -1.

We assume for the rest of this section that $G$ is $p$-valued, and begin the proof of Prop. 6.2. Let $\omega$ be any $p$-valuation on $G$ and define, for each $\nu$ in $\mathbb{R}$, the closed normal subgroups

$$G_{\omega,\nu} = \{g \in G : \omega(g) \geq -\nu\}, \quad G_{\omega,\nu^+} = \{g \in G : \omega(g) > -\nu\}.$$

Put

$$\mathrm{gr}_\omega(G) = \bigoplus_{\nu \in \mathbb{R}} G_{\omega,\nu}/G_{\omega,\nu^+}.$$

The commutator induces a Lie bracket on $\mathrm{gr}_\omega(G)$, which we denote by $[\ ,\ ]_\omega$, and this makes $\mathrm{gr}_\omega(G)$ into a graded Lie algebra over $\mathbb{F}_p$. Let $\mathbb{F}_p[\pi]$ denote the polynomial ring in one variable $\pi$ over $\mathbb{F}_p$, viewed as a graded $\mathbb{F}_p$-algebra with $\pi$ of degree $-1$. The rule $gG_{\nu^+} \mapsto g^p G_{(\nu-1)^+}$ defines an $\mathbb{F}_p$-linear operator $\Pi$ on $\mathrm{gr}_\omega(G)$, which is homogeneous of degree $-1$, and which satisfies $[\Pi g, h]_\omega = \Pi([g, h]_\omega)$ for homogeneous elements $g, h$ of $\mathrm{gr}_\omega(G)$. Letting $\pi$ act as $\Pi$ therefore makes $\mathrm{gr}_\omega(G)$ into a graded Lie algebra over $\mathbb{F}_p[\pi]$.

**Lemma 6.4.** *There exists a $p$-valuation $\omega'$ on $G$ such that, for all $g \in G$ with $g \neq 1$, $\omega'(g) \in e^{-1}\mathbb{Z}$, for some fixed integer $e \geq 1$, and $\mathrm{gr}_{\omega'}(G)$ is an abelian Lie algebra over $\mathbb{F}_p[\pi]$.*

*Proof.* We first convince ourselves that $G$ is of finite rank in the sense of [Laz, Chap. III, Def. 2.1.3], i.e., that $\mathrm{gr}_\omega(G)$ is finitely generated as an $\mathbb{F}_p[\pi]$-module. By [Laz, Chap. III, Prop. 3.1.3 and Prop. 3.1.9] there is an open subgroup $H \subseteq G$ which is $p$-valued of finite rank with respect to the induced $p$-valuation $\omega|H$. Since $\mathrm{gr}_\omega(H)$ is of finite index in $\mathrm{gr}_\omega(G)$ we see that the latter is finitely generated over $\mathbb{F}_p[\pi]$ as well.

In this situation [Laz, Chap. III, Prop. 3.1.11 and Prop. 3.1.12] ensure the existence of a $p$-valuation $\omega'$ on $G$ such that $\omega'(g) \in \mathbb{Q}$ for all $g \neq 1$ and such that $\mathrm{gr}_{\omega'}(G)$ is an abelian Lie algebra. Using [Laz, Chap. III, Prop. 3.1.9] again we obtain that $\mathrm{gr}_{\omega'}(G)$ also is finitely generated over $\mathbb{F}_p[\pi]$. It now follows from [Laz, Chap. III, Prop. 2.2.5 and Prop. 2.2.6] that $\omega'(G \setminus \{1\}) \subseteq a_1 + \mathbb{N}_0 \cup \ldots \cup a_r + \mathbb{N}_0$ for finitely many appropriate rational numbers $a_1, \ldots, a_r$. We finally let $e$ be a common denominator for $a_1, \ldots, a_r$. □

Returning to the proof of Prop. 6.2, we see that, following [Laz, Chap. III, §2.3], our $p$-valuation $\omega'$ on $G$ determines a complete filtration on $\Lambda(G)$, which a priori, is indexed by the additive group $\mathbb{R}$. On the other hand, by [Laz, Chap. III, Thm. 2.3.3], the associated graded ring for $\Lambda(G)$ with this filtration is isomorphic as a graded algebra to the universal enveloping algebra of the $\mathbb{F}_p[\pi]$-Lie algebra $\mathrm{gr}_{\omega'}(G)$. By Lemma 6.4, $\mathrm{gr}_{\omega'}(G)$ is an abelian Lie algebra, which is free over $\mathbb{F}_p[\pi]$ of rank $d = \dim(G)$, and so its universal enveloping algebra is a polynomial ring in $d$ variables over $\mathbb{F}_p[\pi]$. Further, by Lemma 6.4, the degrees of the generators of the universal enveloping algebra belong to $e^{-1}\mathbb{Z}$, and thus Lazard's theorem implies that the filtration on $\Lambda(G)$ is also indexed by $e^{-1}\mathbb{Z}$. Clearly, we can now rescale the filtration on $\Lambda(G)$ so that it is indexed by $\mathbb{Z}$, and the proof of Prop. 6.2 is complete. □

It is now clear that the structure theory as developed in sections 1-5 is applicable to finitely generated modules over the Iwasawa algebra $\Lambda(G)$ of a $p$-valued compact $p$-adic Lie group $G$. The most important open problem remains the classification of the prime c-ideals in $\Lambda(G)$. One might also ask under which additional conditions, $\Lambda(G)$ is a unique factorization ring.

## 7. Examples from elliptic curves

We believe that the abstract theory developed earlier in this paper can be fruitfully applied to many important modules over the Iwasawa algebras of $p$-adic Lie groups, which arise naturally in arithmetic geometry. For brevity, we shall only discuss here one of the most interesting classes of non-commutative examples, which are defined using elliptic curves without complex multiplication. These rather mysterious examples have already been studied (see [CH], [Ho2], [OV]), but it seems certain that a deeper understanding of them, and especially their connection with special values of $L$-functions, will only be achieved by analysing further their structure as modules over the Iwasawa algebra.

Let $F$ be a finite extension of $\mathbb{Q}$, and $E$ an elliptic curve defined over $F$. We shall always assume that $E$ has no complex multiplication (i.e. that the endomorphism ring of $E$ over $\overline{\mathbb{Q}}$ is $\mathbb{Z}$). Let $p$ be any prime number $\geq 5$. For each integer $n \geq 1$, we write $E_{p^n}$ for the group of $p^n$-division points on $E$, and $E_{p^\infty}$ for the group of all $p$-power division points on $E$. We define

$$F_\infty = F(E_{p^\infty}), \quad G = G(F_\infty/F).$$

The action of $G$ on $E_{p^\infty}$ defines an injection of $G$ into $\mathrm{Aut}(E_{p^\infty}) \cong GL_2(\mathbb{Z}_p)$; and a celebrated theorem of Serre [Se2] asserts that $G$ is open in $GL_2(\mathbb{Z}_p)$. Let $\mu_{p^\infty}$ denote the subgroup of all $p$-power roots of unity. By the Weil pairing, the field $F(\mu_{p^\infty})$ is contained in $F_\infty$, and we write $\Gamma$ for the Galois group of $F(\mu_{p^\infty})$ over $F$. Henceforth, we shall always assume that the Galois group $G$ is pro-$p$ (this can always be achieved by replacing $F$ by a finite extension, if necessary). This implies that $\Gamma$ is pro-$p$, and hence that $F(\mu_{p^\infty})$ is itself the cyclotomic $\mathbb{Z}_p$-extension of $F$. For simplicity, we write $F^{\mathrm{cyc}} = F(\mu_{p^\infty})$.

The modules which interest us arise as follows from the classical Selmer groups, whose definition we now recall. Let $L$ denote any intermediate field with $F \subseteq L \subseteq F_\infty$. As usual, if $v$ is any finite place of $L$, we write $L_v$ for the union of the completions at $v$ of all finite extensions of $F$ contained in $L$. If $K$ is a field, $\overline{K}$ will always denote a fixed algebraic closure of $K$.

**Definition 7.1.** The *Selmer group* $\mathcal{S}(E/L)$ is defined as

$$\mathcal{S}(E/L) = \mathrm{Ker}(H^1(G(\overline{\mathbb{Q}}/L), E_{p^\infty}) \to \prod_v H^1(G(\overline{L_v}/L_v), E(\overline{L_v}))).$$

We also write $X(E/L)$ for the compact Pontryagin dual

$$X(E/L) = \mathrm{Hom}(\mathcal{S}(E/L), \mathbb{Q}_p/\mathbb{Z}_p).$$

If $L$ is Galois over $F$, then $G(L/F)$ has a natural action on both $\mathcal{S}(E/L)$ and $X(E/L)$, and this extends to a natural action of $\Lambda(G(L/F))$ of the $p$-adic Lie group $G(L/F)$. Moreover, it is easy to see that $X(E/L)$ is a finitely generated module over $\Lambda(G(L/F))$. Lying much deeper, and still only proven in a very modest number of cases, are the following two conjectures. If $E$ has good ordinary reduction at all primes $v$ of $F$ dividing $p$, it is conjectured that $X(E/F^{\mathrm{cyc}})$ is always $\Lambda(\Gamma)$-torsion (Mazur [Ma]), and that $X(E/F_\infty)$ is $\Lambda(G)$-torsion (Harris [Ha], [CH]). We mention the two conjectures at the same time because one can study both by exploiting the connection between them. Put

$$H = G(F_\infty/F^{\mathrm{cyc}}),$$

so that $\Lambda(H)$ is a sub-algebra of $\Lambda(G)$. Then $\Lambda(G)$ is not finitely generated as a $\Lambda(H)$-module, because $\Gamma = G/H$ is infinite. Now the structure theory of finitely generated torsion $\Lambda(\Gamma)$-modules is very well-known (see [B-CA], [Iw]). Indeed, if $M$ is a finitely generated $\Lambda(\Gamma)$-module, this structure theory shows that $M$ is a finitely generated $\mathbb{Z}_p$-module if and only if $M$ is $\Lambda(\Gamma)$-torsion and the $\mu$-invariant of $M$ is zero. The following result is proven in [CH].

**Proposition 7.2.** *Assume the following hypotheses are valid:*
*(i) $p \geq 5$,*
*(ii) $G$ is pro-$p$,*
*(iii) $E$ has good ordinary reduction at all places $v$ of $F$ dividing $p$, and*
*(iv) $X(E/F^{\text{cyc}})$ is a finitely generated $\mathbb{Z}_p$-module.*
*Then $X(E/F_\infty)$ is a finitely generated $\Lambda(H)$-module, and so, in particular, $X(E/F_\infty)$ is $\Lambda(G)$-torsion.* □

We mention in passing that it is not always true that $X(E/F_\infty)$ is finitely generated over $\Lambda(H)$ (many examples are known, see [CH], where $X(E/F^{\text{cyc}})$ has positive $\mu$-invariant, and for these examples $X(E/F_\infty)$ cannot be finitely generated over $\Lambda(H)$).

A second basic result about $X(E/F_\infty)$ is due to Ochi and Venjakob [OV].

**Proposition 7.3.** *Assume that hypothesis (i), (ii) and (iii) of Prop. 7.2 are valid. If $X(E/F_\infty)$ is $\Lambda(G)$-torsion, then it contains no non-zero pseudo-null submodule. In particular, if $X(E/F_\infty)$ is finitely generated over $\Lambda(H)$, then $X(E/F_\infty)$ has no non-zero pseudo-null submodule.* □

**Corollary 7.4.** *Assume that the hypotheses (i), (ii), (iii) and (iv) of Prop. 7.2 are valid. Then the $\Lambda(H)$-torsion submodule of $X(E/F_\infty)$ is zero, and $X(E/F_\infty)$ has positive $\Lambda(H)$-rank.*

*Proof.* It is proven in [OV] that every $\Lambda(G)$-module which is finitely generated and torsion over $\Lambda(H)$ is pseudo-null as a $\Lambda(G)$-module. For the final assertion, we only need to note that it is shown in [CH] that $X(E/F_\infty) \neq 0$ for all primes $p \geq 5$. □

We cannot resist mentioning a curious arithmetic consequence of Cor. 7.4, whose interest lies in the fact that it asserts a non-triviality result for all good ordinary primes $p \geq 5$. For this proposition alone, we drop our assumption that $G$ is pro-$p$.

**Proposition 7.5.** *Let $E$ be an elliptic curve over $F$ without complex multiplication, whose j-invariant is an algebraic integer. Let $p$ be any prime number $\geq 5$ such that $E$ has good ordinary reduction at all primes $v$ of $F$ dividing $p$. Put $K = F(E_p)$. Then $\mathcal{S}(E/K^{\text{cyc}})$ is infinite.*

*Proof.* Since the j-invariant of $E$ is an integer in $F$, $E$ has good reduction everywhere over $K = F(E_p)$ by the results of [ST]. Put $G_K = G(F_\infty/K)$, so that $G_K$ is pro-$p$. Let us assume that $\mathcal{S}(E/K^{\text{cyc}})$ is finite, and derive a contradiction. Now Cor. 7.4 implies that $X(E/F_\infty)$ is a finitely generated $\Lambda(H_K)$-module, of strictly positive $\Lambda(H_K)$-rank, which we denote by $r$; here $H_K = G(F_\infty/F^{\text{cyc}})$. Now let $L$ be a variable finite extension of $K$ contained in $F_\infty$, and put $H_L = G(F_\infty/L^{\text{cyc}})$. A well-known algebraic argument (see [Ho0]) then shows that, as $[L^{\text{cyc}} : K^{\text{cyc}}] \to \infty$, we have the asymptotic estimate

$$\mathbb{Z}_p\text{--rank of } X(E/F_\infty)_{H_L} = r[L^{\text{cyc}} : K^{\text{cyc}}] + o([L^{\text{cyc}} : K^{\text{cyc}}]).$$

Here $X(E/F_\infty)_{H_L}$ denotes the $H_L$-coinvariants of $X(E/F_\infty)$. We shall only need the very weak consequence of this result that $X(E/F_\infty)_{H_L}$ is not finite when $[L^{\text{cyc}} : K^{\text{cyc}}]$ is sufficiently large. On the other hand, if $L$ is any finite Galois extension of $K$ contained in $F_\infty$, we can apply the formula of Hachimori and Matsuno (see [HM], or [Ho3, Cor. 2.12]) to the finite Galois $p$-extension $L^{\text{cyc}}/K^{\text{cyc}}$ to conclude, on recalling that $E$ has good reduction everywhere over $K$, that $\mathcal{S}(E/L^{\text{cyc}})$ is finite always. But the restriction map from $\mathcal{S}(E/L^{\text{cyc}})$ to $\mathcal{S}(E/F_\infty)^{H_L}$ induces a homomorphism

$$f : X(E/F_\infty)_{H_L} \to X(E/L^{\text{cyc}}).$$

A non-trivial arithmetic argument (see [CH, §6, Lemma 6.7]) shows that $\text{Ker}(f)$ and $\text{Coker}(f)$ are both finite, since $E$ has no places of bad multiplicative reduction. Thus we conclude that $X(E/F_\infty)_{H_L}$ is finite for all finite Galois extensions $L$ of $K$ contained in $F_\infty$, contradicting what we proved above. This completes the proof of the proposition. □

**Example 7.6.** Let $E$ be the elliptic curve over $\mathbb{Q}$ given by

$$E : y^2 + xy = x^3 + x^2 - 2x - 7.$$

The conductor of $E$ is 121 ($E$ is the curve 121C1 in [Cr]), and the j-invariant of $E$ is equal to $-11^2$. We mention in passing that Serre [Se2] has proven that $G$ is isomorphic to $GL_2(\mathbb{Z}_p)$ for all primes $p \neq 11$, where $G = G(\mathbb{Q}(E_{p^\infty})/\mathbb{Q})$. The good ordinary primes for $E$ are $2, 3, 5, 7, 13, 17, \cdots$ (the first supersingular prime is 43). Thus Prop. 7.5 can be applied to $E$ and proves that for all ordinary primes $p \geq 5$, the Selmer group $\mathcal{S}(E/\mathbb{Q}(E_p, \mu_{p^\infty}))$ is infinite. Put $F = \mathbb{Q}(E_p)$. To illustrate the limits of our present knowledge, we still cannot prove that either $X(E/F^{\text{cyc}})$ is $\Lambda(\Gamma)$-torsion or $X(E/F_\infty)$ is $\Lambda(G)$ -torsion for a single good ordinary prime $p \geq 5$ for this curve $E$.

It is now natural to try and compare these arithmetic results about the module $X(E/F_\infty)$ with the algebraic theory established earlier in this paper. All we can do at present is to pose the following rather obvious questions. For simplicity, let us assume that our elliptic curve $E$ over $F$ and the prime $p$ satisfy hypotheses (i) - (iv) of Prop. 7.2, so that we know that $X(E/F_\infty)$ is $\Lambda(G)$-torsion. In addition, let us assume that $G$ is a $p$-valued group in the sense of §6 (if this is not the case, it can always be achieved by replacing $F$ by a finite extension contained also in $F_\infty$). Then, by the results of §6, $A = \Lambda(G)$ is a noetherian Auslander regular maximal order, without zero divisors, satisfying $\text{K}_1\dim(A) \leq 1$. Take $M = X(E/F_\infty)$, and as before, we write $Q(M)$ for the corresponding object in the quotient category $\text{Mod}(A)/\mathcal{C}_A^1$, where $A = \Lambda(G)$. By Prop. 4.1, we have the canonical decomposition

$$Q(M) = Q(M)_0 \oplus Q(M)_1 ,$$

where $Q(M)_0$ is completely faithful, and $Q(M)_1$ is bounded. At present, we do not know a single example of an elliptic curve $E$ over $F$ and a prime $p$ satisfying our hypothesis for which we can prove that a given one of the direct summands $Q(M)_i$ ($i = 0, 1$) is non-zero. Of course, one at least of the two direct summands must be non-zero, since Prop. 7.3 and Cor. 7.4 show that $Q(M)$ is non-zero. It follows from Cor. 2.5 and Prop. 7.3 that there exist non-zero left ideals $L_1, \cdots, L_m$ of $\Lambda(G)$ such that we have an exact sequence of $\Lambda(G)$-modules

$$0 \to \bigoplus_{i=1}^m \Lambda(G)/L_i \to X(E/F_\infty) \to Z \to 0,$$

where $Z$ is pseudo-null. We stress that there appears to be great interest arithmetically in studying the left ideals $L_1, \cdots, L_m$. One can hope eventually to give an analytic description of the characteristic ideal of $X(E/F_\infty)$ in terms of the values at $s = 1$ of the twists of the complex $L$-function of $E$ over $F$ by Artin characters of $G$ (i.e. those characters of $G$ which factor through a finite quotient). A more modest goal at present is to attempt to determine these ideals in some simple numerical examples We now discuss one such example in more detail.

**Example 7.7.** Let $E$ be the elliptic curve $X_1(11)$, namely $E : y^2 + y = x^3 - x^2$, of conductor 11. Take $p = 5$, and let $F$ be the field obtained by adjoining the 5-th roots of unity to $\mathbb{Q}$. Thus $F_\infty = \mathbb{Q}(E_{5^\infty})$, and $F^{\text{cyc}} = \mathbb{Q}(\mu_{5^\infty})$, and we have

$$G = G(F_\infty/F), \quad H = G(F_\infty/\mathbb{Q}(\mu_{5^\infty})).$$

It can easily be shown (see [Fi]) that the image of $G$ in $\mathrm{Aut}(E_{5^\infty})$ can be identified with the subgroup of all matrices $\begin{pmatrix} a & b \\ c & d \end{pmatrix}$ in $GL_2(\mathbb{Z}_5)$ with $a \equiv d \equiv 1 \mod 5$, and $c \equiv 0 \mod 5^2$, and this group in turn is isomorphic to the group of all matrices in $GL_2(\mathbb{Z}_5)$ which are congruent to the identity modulo 5. Hence $G$ is pro-5, and also a 5-valued group in the sense of §6. As 5 is an ordinary prime for $E$, and it is well-known that $\mathcal{S}(E/\mathbb{Q}(\mu_{5^\infty})) = 0$ (see [CS, Chap. 5]), we see that hypotheses (i)-(iv) of Prop. 7.2 are valid in this case.

**Proposition 7.8** *In the above example, $X(E/F_\infty)$ is a torsion $\Lambda(G)$−module, with no non-zero pseudo-null submodule. It is finitely generated over $\Lambda(H)$ of rank 4, its $\Lambda(H)$-torsion submodule is zero, but it is not a free $\Lambda(H)$-module.* □

The proof that $X(E/F_\infty)$ has $\Lambda(H)$-rank 4 hinges on the remarkable fact that we can determine the exact $\mathbb{Z}_5$-rank of the $H_L$-coinvariants of $X(E/F_\infty)$ for any finite Galois extension $L$ of $F$ contained in $F_\infty$, where $H_L = G(F_\infty/L^{\mathrm{cyc}})$. The following result which depends crucially on the ideas of [HM], summarizes what is proven in §7 of [CH].

**Lemma 7.9** *Let $L$ be any finite Galois extension of $F$ contained in $F_\infty$. Then*
*(i) $X(E/L^{\mathrm{cyc}})$ is a free $\mathbb{Z}_5$-module of rank $4 \cdot [L^{\mathrm{cyc}} : F^{\mathrm{cyc}}] - r_L$, where $r_L$ denotes the number of primes of $L^{\mathrm{cyc}}$ above 11.*
*(ii) $X(E/F_\infty)_{H_L}$ has $\mathbb{Z}_5$-rank equal to $4 \cdot [L^{\mathrm{cyc}} : F^{\mathrm{cyc}}]$, and its $\mathbb{Z}_5$-torsion submodule contains a group having the same order as $E_{5^\infty}(L^{\mathrm{cyc}})$.* □

The only new ingredient in the proof of Lemma 7.9 is the result of Matsuno [Mat] that $X(E/L^{\mathrm{cyc}})$ contains no $\mathbb{Z}_5$-torsion, and we do not repeat the proof here. Our hope is that Lemma 7.9 might eventually enable us to say more about the questions on the structure of $X(E/F_\infty)$ raised earlier. However, it seems likely that one will have to establish some result about the action of the centre of $G$ on $X(E/F_\infty)$ and this does not appear to be easy.

By contrast to the above calculation of the exact $\mathbb{Z}_5$-rank of $X(E/F_\infty)_{H_L}$, it seems to be a deep arithmetic problem to determine the exact $\mathbb{Z}_5$-rank of $Y(L) := X(E/F_\infty)_{G_L}$, where again $L$ ranges over the finite extensions of $F$ contained in $F_\infty$, and $G_L = G(F_\infty/L)$. This is hardly surprising since the methods of [CH] show that the $\mathbb{Z}_5$-rank of $Y(L)$ is equal to the $\mathbb{Z}_5$-rank of $X(E/L)$. In particular, $Y(L)$ is finite if and only if both $E(L)$ and the 5-primary part of the Tate-Shafarevich group of $E$ over $L$ are finite. As remarked above, we have already used the classical fact (see [CS, Chap. 5]) that $Y(F)$ is finite. By some remarkable explicit descent calculations, Fisher [Fi] has recently shown that $Y(L)$ is finite for $L$ ranging over the following four Galois extensions of $\mathbb{Q}$, which are cyclic of degree 5 over $F$ and contained in $F_\infty$,

$$L_1 = \mathbb{Q}(X_1(11)_5) , \quad L_2 = \mathbb{Q}(X_2(11)_5) , \quad L_3 = F\mathbb{Q}(\mu_{11})^+ , \quad L_4 = \mathbb{Q}(\mu_{25}) ;$$

here $X_2(11)$ denotes the unique elliptic curve of conductor 11 over $\mathbb{Q}$ which has no non-zero rational point, and $\mathbb{Q}(\mu_{11})^+$ denotes the maximal real subfield of the field of 11-th roots of unity. Moreover, the first author, Fisher, and Greenberg have deduced from Fisher's results by arguments from Iwasawa theory that $Y(L)$ is finite for every finite extension $L$ of $F$ contained in the field obtained by adjoining to $L_1$ all 5-power roots of unity. These fragmentary results lead one to raise the mysterious question of whether $Y(L) = X(E/F_\infty)_{G_L}$ is finite for all finite extensions $L$ of $F$ contained in $F_\infty$.

We end this paper by pointing out that there is a possible simple explanation, why the above question should have a positive answer, in terms of the structure theory, thanks to the following general result. Let $G$ be a compact $p$-adic Lie group which is $p$-valued. By an *Artin representation* of $G$, we mean a homomorphism

$$\rho : G \longrightarrow \mathrm{Aut}(W) ,$$

where $W$ is a finite dimensional vector space over $\bar{\mathbb{Q}}_p$, which satisfies $\rho(U) = 1$ for some open normal subgroup $U$ of $G$. Let $\rho$ be such an Artin representation of $G$. By the universal property of the group ring $\Lambda(G/U)$ it extends to a homomorphism of rings

$$\rho : \Lambda(G) \longrightarrow \Lambda(G/U) \longrightarrow \mathrm{End}(W) .$$

We can then define the map

$$\det(\rho) := \det \circ \rho : \Lambda(G) \longrightarrow \bar{\mathbb{Q}}_p .$$

**Lemma 7.10** *Assume that $M$ is a finitely generated $\Lambda(G)$-module, which has no non-zero pseudo-null submodule. Assume further that there exists an element $a \in \mathrm{ann}(Q(M))$ such that $\det(\rho)(a) \neq 0$ for all Artin representations $\rho$ of $G$. Then for all open subgroups $U \subseteq G$ the $U$-coinvariants $M_U$ are finite.*

*Proof.* Since $M$ has no non-zero pseudo-null submodule, it follows from [Rob, Lemma 2.5] that $\mathrm{ann}(Q(M)) = \mathrm{ann}_{\Lambda(G)}(M)$. By replacing $U$ by a smaller subgroup if necessary, we can assume that $U$ is also normal in $G$. As the element $a$ belongs to $\mathrm{ann}_{\Lambda(G)}(M)$, it is clear that the image $a_U$ of $a$ in $\Lambda(G/U)$ annihilates $M_U$. On the other hand, $M_U$ is a finitely generated $\Lambda(G/U)$-module because $M$ is a finitely generated $\Lambda(G)$-module. Hence $M_U$ is a finitely generated $\Lambda(G/U)/\Lambda(G/U)a_U$-module. But this latter ring $\Lambda(G/U)/\Lambda(G/U)a_U$ is finite because $\det(\rho)(a) \neq 0$ for all irreducible Artin representations $\rho$ of $G$ with $\rho(U) = 1$. This completes the proof of the lemma. $\square$

John Coates
DPMMS, University of Cambridge
Centre for Mathematical Sciences
Wilberforce Road
Cambridge CB3 0WB, England
j.h.coates@pmms.cam.ac.uk

Peter Schneider
Mathematisches Institut
Westfälische Wilhelms-Universität Münster
Einsteinstr. 62
D-48149 Münster, Germany
pschnei@math.uni-muenster.de
http://www.uni-muenster.de/math/u/schneider

R. Sujatha
School of Mathematics
Tata Institute of Fundamental Research
Homi Bhabha Road,
Colaba
Mumbai 400005
India
sujatha@math.tifr.res.in